\documentclass[12pt]{amsart}
\usepackage{amsfonts, amscd, amsmath, amssymb}

\newcommand\globcnt{subsubsection}

\swapnumbers
\theoremstyle{plain}

\newtheorem{thm}[\globcnt]{Theorem}
\newtheorem{lem}[\globcnt]{Lemma}
\newtheorem{prop}[\globcnt]{Proposition}
\newtheorem{cor}[\globcnt]{Corollary}

\newtheorem{claim}[\globcnt]{Claim}
\theoremstyle{definition}
\newtheorem{defn}[\globcnt]{Definition}

\newcommand\MainTheorem{Theorem~\ref{th:MainTheorem}}

\newcommand\partA{(A)}
\newcommand\partB{(B)}

\newcommand\proofstyle[1]{{\it #1}}
\providecommand\proof{\proofstyle{Proof.}\ }

\newcommand\tangst[2]{T{#1}_{#2}}
\newcommand\tanpm{\tangst{\manif}{\pnt}}    

\newcommand\CCC{{\mathbb C}}
\newcommand\NNN{{\mathbb N}}
\newcommand\RRR{{\mathbb R}}

\newcommand\ZZZ{{\mathbb Z}}

\newcommand\JJ{{\mathcal J}}

\newcommand\MM{{\mathcal M}}
\newcommand\NN{{\mathcal N}}

\newcommand\UU{{\mathcal U}}

\newcommand\XX{{\mathcal X}}

\newcommand\NNNz{\NNN_{0}}
\newcommand\NNNii{\overline{\NNN}}
\newcommand\NNNzi{\NNNii_{0}}

\newcommand\Int{\mathrm{Int}}

\newcommand\id{\mathrm{id}}

\newcommand\Fr{\mathrm{Fr}}

\newcommand\Per{\mathrm{Per}}
\newcommand\IM{\mathrm{Im}}

\newcommand\diag{\mathrm{diag}}

\newcommand\vect[2]{\left\|\begin{array}{c}#1\\#2\end{array}\right\|}

\newcommand\matr[4]{\left\|\begin{array}{cc}#1 & #2 \\ #3 & #4 \end{array}\right\|}

\newcommand\eps{\varepsilon}

\newcommand\unit[1]{{#1}_{0}}                 

\newcommand\surf{M} 
\newcommand\manif{\surf}

\newcommand\Fix{{\mathrm{Fix\,}}}

\newcommand\Diff{\mathcal{D}}             
\newcommand\DiffId{\unit{\Diff}}
\newcommand\DiffM{\Diff(\manif)}

\newcommand\Fld{F}

\newcommand\flow{\Phi}
\newcommand\tflow{\widetilde{\flow}}
\newcommand\Shflow{\Omega}

\newcommand\Frechet{Fr\'echet}  

\newcommand\ImFn{\Gamma_{\flow}}
\newcommand\ImFnplus{\ImFn^{+}}
\newcommand\ImFnminus{\ImFn^{-}}

\newcommand\jet[3]{\left[#2\right]^{#1}_{#3}}

\newcommand\NbhW[3]{\NN_{\,#2}^{\,#3}(#1)}

\newcommand\NbhS[2]{\NN^{\,#2}(#1)}

\newcommand\Rot[1]{R(#1)}

\newcommand\smsp[1]{C^{\infty}\ifx#1\empty\else(#1)\fi}
\newcommand\smfunc{\smsp{\manif,\RRR}}
\newcommand\smmap{\smsp{\manif, \manif}}

\newcommand\afunc{\alpha}
\newcommand\tafunc{\tilde\alpha}
\newcommand\Afunc{\Lambda}
\newcommand\tAfunc{\tilde\Afunc}

\newcommand\amap{f}
\newcommand\tamap{\tilde\amap}
\newcommand\Amap{F}

\newcommand\bfunc{\beta}
\newcommand\bmap{g}

\newcommand\cfunc{\gamma}
\newcommand\cmap{h}

\newcommand\Dflow{\Diff(\flow)}
\newcommand\DflowId{\DiffId(\flow)}

\newcommand\End{{\mathcal{E}}}
\newcommand\EndId{\unit{\End}}

\newcommand\Cflow{\End(\flow)}
\newcommand\CflowId{\EndId(\flow)}

\newcommand\condExt{(E)}
\newcommand\condExtk{(E)${}^{k}$}
\newcommand\condExtzero{(E)${}^{0}$}
\newcommand\condExtnk{(E)${}^{n+k}$}

\newcommand\B{B}

\newcommand\tUU{\widetilde\UU}
\newcommand\UUi{\UU_{i}}
\newcommand\tUUi{\tUU_{i}}

\newcommand\pnt{z}

\newcommand\pntb{y}
\newcommand\pntc{z}

\providecommand\nb{U}
\providecommand\nbs{V}
\providecommand\nbi{U_i}

\newcommand\nbb{V}
\newcommand\nbc{W}

\newcommand\NBm{\NN}
\newcommand\NBf{\MM}

\newcommand\unifunc{\delta}
\newcommand\unifuncA{\widehat{\unifunc}}

\newcommand\NBfd{\NBf^{\unifunc}}
\newcommand\NBfdA{\NBf^{\unifuncA}}

\newcommand\NBfnbi{\NBf_{i}}
\newcommand\tNBfnbi{\NBf^{\,\prime}_{i}}

\newcommand\NBfdnbs{\NBf^{\unifunc}_{\nbs}}

\newcommand\NBmi{\NBm_{i}}
\newcommand\tNBmi{{\widetilde\NBm}_{i}}
\newcommand\tNBm{\widetilde\NBm}

\newcommand\MapOnSet[2]{{#1}_{#2}}              
\newcommand\ShM[1]{\MapOnSet{\Shift}{#1}}

\newcommand\condS{(S)}
\newcommand\condSk{(S)${}^{k}$}
\newcommand\condSz{(S)${}^{0}$}
\newcommand\condSnk{(S)${}^{n+k}$}
\newcommand\rMmapi{p_i}

\newcommand\smk{r}
\newcommand\epsi{\eps_{i}}
\newcommand\Ki{K_{i}}

\newcommand\Z{Z}
\newcommand\Zid{\Z_{\id}}

\newcommand\rcomp{C}
\newcommand\nbhcomp{V}
\newcommand\nind{\lambda}
\newcommand\im{\mathrm{Im}}
\newcommand\GL{GL}

\newcommand\Odzflow{U}
\newcommand\Interval{\JJ}
\newcommand\Shift{\varphi}
\newcommand\smcl{r}

\newcommand\sectflow{\smsp{\Odzflow,\Interval}}
\newcommand\shiftflow{\smsp{\Odzflow,\manif}}

\newcommand\inb{\id_{\Odzflow}}
\newcommand\nbhtmp{W}

\newcommand\extflow{\Phi}
\newcommand\factflow{\Psi}
\newcommand\extflsh{\phi}
\newcommand\factflsh{\psi}
\newcommand\tp{P}

\newcommand\Diskk{D^{k}}

\newcommand\mult{Z}
\newcommand\norm[3]{{\left\| #1 \right\|}_{#2,#3}}
\newcommand\normnK[1]{\norm#1{\smk}{K}}
\newcommand\normnprevK[1]{\norm#1{\smk-1}{K}}
\newcommand\diam{\mathrm{diam}}

\newcommand\stmap{\cmap}
\newcommand\stfunc{\sigma}

\newcommand\dmap{\cfunc}
\newcommand\onb{\overset{\circ}{\nb}}
\newcommand\Field{\mathbb{F}}

\newcommand\Bmap{H}

\newcommand\DimNbhnb[1]{{\nb}^{#1}}
\newcommand\nbm{\DimNbhnb{m}}
\newcommand\nbn{\DimNbhnb{n}}
\newcommand\nbk{\DimNbhnb{k}}
\newcommand\nbmn{\DimNbhnb{m+n}}
\newcommand\nbnk{\DimNbhnb{n+k}}
\newcommand\nbmnk{\DimNbhnb{m+n+k}}

\newcommand\unifuncext{\unifunc_{\extflsh}}
\newcommand\unifuncfact{\unifunc_{\factflsh}}

\newcommand\NBffact{\NBf_{\factflsh}}
\newcommand\NBfext {\NBf_{\extflsh}}

\newcommand\tNN{\widetilde\NN}

\newcommand\dflow{\Psi}
\newcommand\expflow{\nu}

\newcommand\teyl{\nu} 

\newcommand\barz{\bar z}
\newcommand\baromega{\bar\omega}

\newcommand\I{ {\mathcal I} }

\newcommand\dimM{m}

\newcommand\Inv{\mathrm{Inv}}
\newcommand\Kinv{\Inv}
\newcommand\FlowInv{\Kinv(\flow)}
\newcommand\FlowInvId{\Kinv_{0}(\flow)}

\newcommand\numb{A}
\newcommand\tm{t}

\begin{document}

\title{Smooth shifts along trajectories of flows}
\author{Sergey Maksymenko}
\email{maks@imath.kiev.ua}
\address{
Topology Division,
Institute of Mathematics of Ukrainian NAS,
Tereshchenkovskaya str. 3,
01601 Kiev,
Ukraine
}

\begin{abstract}
Let $\flow$ be a flow on a smooth, compact, finite-dimensional manifold $\manif$.
Consider the subset $\Dflow$ of $\smmap$
consisting of diffeomorphisms of $\manif$ preserving the foliation of the flow $\flow$.
Let also $\DflowId$ be the identity path component of $\Dflow$ with compact-open topology.
We prove that under mild conditions on fixed points of $\flow$
the space $\DflowId$ is either contractible or homotopically equivalent to $S^1$.
\end{abstract}



\maketitle

\section{Introduction}
Throughout, $\manif$ will be a smooth ($C^{\infty}$) connected manifold
and $\DiffM$ be the space of diffeomorphisms of $\manif$.
Let $\flow$ be a smooth flow on $\manif$, and
$\Fix\flow$ be the fixed-point set of $\flow$.
Define the map $\Shift:\smfunc\to\smmap$ by
$$
 \Shift(\afunc)(x) = \flow(x, \afunc(x)),
$$
where $\afunc\in\smfunc$ and $x\in\manif$.
We will say that $\Shift$ is the {\em shift-map\/} along trajectories of $\flow$.
If $\afunc\in\smfunc$ and $\amap = \Shift(\afunc)$,
then the following statements can easily be checked:
\begin{enumerate}
\item
 $\amap$ is homotopic to $\id_{\manif}$;
\item
 $\amap(\omega)\subset\omega$
 for each trajectory $\omega$ of $\flow$. In particular,
 if $\pnt\in\Fix\flow$, then $\amap(\pnt)=\pnt$. Moreover,
\item
 $\amap$ is a local diffeomorphism at $\pnt$, i.e.,
 the corresponding tangent map $\amap'(\pnt):\tanpm\to \tanpm$ is an isomorphism.
\end{enumerate}
Let $\Cflow \subset \smmap$ be the set of all maps $\amap:\manif\to\manif$
satisfying the above conditions (2) and (3),
$\Dflow$ be the intersection $\Cflow \cap \DiffM$,
and $\CflowId$ and $\DflowId$ be the identity path components of the spaces
$\Cflow$ and $\Dflow$ (respectively) in the compact-open topology.
Evidently $\IM\Shift \subset \CflowId$.

We use the map $\Shift$ to study the homotopy types of
the spaces $\CflowId$ and $\DflowId$.
Take any $\smcl=0,1,\ldots,\infty$ and endow $\smfunc$ and $\smmap$ with
the strong Whitney $C^{\smcl}$-topologies and
the spaces $\IM\Shift$, $\DflowId$, and $\CflowId \subset \smmap$
with the induced ones.
The following theorem summarizes the results obtained in the paper.

\begin{thm}\label{th:MainTheorem}
Suppose for each fixed point $\pnt$ of $\flow$ there exist local coordinates
$(x_1,\ldots,x_{\dimM})$ and a linear flow $\factflow$ on $\RRR^{n} (n \leq \dimM)$ such
that $\pnt=0\in\RRR^{\dimM}$ and for all $t$ in some neighborhood of $0\in\RRR$,
we have $p \circ \flow_t = \factflow_t \circ p$,
where $p:\RRR^{\dimM} \to \RRR^n$ is the natural projection.
Then
\begin{enumerate}
\item[{\em \partA}]
$\im\Shift = \CflowId$ and $\Shift:\smfunc \to \CflowId$ is
either a homeomorphism or
a covering map with $\ZZZ$ as a group of covering slices.
Also, the set $\Shift^{-1}(\DflowId)$ is {\em convex\/} if regarded as a subset of
the {\em linear\/} space $\smfunc$.
\item[{\em \partB}]
If $\manif$ is compact,
then the inclusion $\DflowId \subset \CflowId$ is a homotopy equivalence
and these spaces are either contractible or have the homotopy type of $S^{1}$.
They are contractible whenever
$\flow$ has at least one non-closed trajectory or
if the tangent linear flow at some fixed point of $\flow$ is trivial.
\end{enumerate}
\end{thm}

There are some applications of this
theorem to Morse functions and Morse-Smale flows
similar to~\cite{Gutierrez_Melo}.
We will consider them in the next papers.

The paper is organized as follows.
In section~\ref{sect:Preliminaries} we recall the definitions
of Whitney topologies and some formulas concerning linear flows.

Section~\ref{sect:map_gen_flow}.
We show that the set $\Zid = \Shift^{-1}(\id_{\manif})$
is a {\em subgroup\/} of $\smfunc$ and for each $\afunc\in\smfunc$, we have
$\Shift^{-1} \Shift(\afunc) = \{ \afunc + \Zid \}$.
Thus the correspondence $\{\afunc+\Zid\} \mapsto \Shift(\afunc)$,
where $\afunc\in\smfunc$,
is a bijection of the factor-group $\smfunc/\Zid$ onto the image $\IM\Shift$.
Notice, however, that $\Shift$ is not a homomorphism
(see formulas~\eqref{equ:sh_gh} and~\eqref{equ:sh_g_1}).
We also describe $\Zid$ in terms of the interior of $\Fix\flow$
(Theorem~\ref{th:Zid_descr}).
In particular, we obtain that $\Fix\flow$ is
nowhere dense in $\manif$, i.e., $\Int\Fix\flow=\emptyset$ if and only if
$\Zid$ is either trivial or isomorphic to $\ZZZ$.

Section~\ref{sect:loc_sect_of_invsh}.
There are two natural topologies on $\IM\Shift$: the factor-topology of $\smfunc$
and the induced topology of the ambient space $\smmap$.
Assuming $\Int\Fix\flow=\emptyset$, we prove that under mild
conditions on fixed points of $\flow$,
these topologies coincide (Theorem~\ref{th:LocSect}).
In this case $\Shift:\smfunc\to\IM\Shift$ is either a homeomorphism
or a covering map.

Section~\ref{sect:diff_gen_flow}.
We study the set $\ImFn = \Shift^{-1}(\DiffM)$.
It is a union of two disjoint, open subsets $\ImFnplus$
and $\ImFnminus$ corresponding
to diffeomorphisms of $\manif$ preserving and reversing (respectively)
orientation of trajectories of $\flow$.
We prove that
$\ImFnplus$ and $\ImFnminus$ are {\em convex\/} if regarded as subsets of
the {\em linear\/} space $\smfunc$ (Theorem~\ref{th:h_diff}).
Also, if $\amap \in \IM\Shift$ and $\pnt\in\Fix\flow$,
then the tangent map $\amap'(\pnt)$ is an isomorphism
whence $\IM\Shift \subset \CflowId$.

Section~\ref{sect:image_of_invsh}.
A sufficient condition for the relation $\IM\Shift = \CflowId$ is
given (Theorem~\ref{th:IMinvsh=CflowId}).
In this case we have $\Shift(\ImFnplus) = \DflowId$.

Sections~\ref{sect:calculations} and~\ref{sect:proof_HadamLemmas}.
We show that Theorems~\ref{th:LocSect} and~\ref{th:IMinvsh=CflowId} hold true
for a flow $\flow$ that satisfies the conditions of \MainTheorem\
(Theorem~\ref{th:shift_reg_ext_lin}).
This proves part \partA\ of \MainTheorem.

Section~\ref{sect:HomotopyGroups}.
Assuming $\manif$ is compact, we describe the homotopy
type of the spaces $\DflowId$ and $\CflowId$.
This completes \MainTheorem.

Finally, in section~\ref{sect:CflowId_closure}
we shortly discuss the closure of $\CflowId$.

The author wishes to express his indebtedness to
V.~V.~Sharko, H.~Zieschang, M.~Pankov, E.~Polulyah, A.~Prishlyak, I.~Vlasenko,
E.~Kud\-ryavt\-zeva, and O.~Mozgova for many helpful conversations.
The author thanks the referee
for referring him to the paper of J.~Keesling
and for pointing to the closures of the sets $\DflowId$ and $\CflowId$.
The author is also grateful to Yu.~A.~Chapovsky for careful reading
this manuscript and correcting numerous misprints.

\section{Preliminaries}\label{sect:Preliminaries}
\subsection{Whitney topologies}\label{sect:WhitneyTopologies}
In the sequel, all manifolds are assumed to be smooth ($C^{\infty}$).
Let us denote
$\NNNz=\NNN\cup\{0\}$ and $\NNNzi=\NNN\cup\{0,\infty\}$.
For a subset $X$ of a topological space $Y$ the symbols
$\overline{X}$, $\Int X$, and $\Fr X = \overline{X}\setminus\Int X$
mean the closure, the interior,
and the frontier of $X$ in $Y$, respectively.

Let $M$ and $N$ be manifolds,
$f \in \smsp{M,N}$, $x \in M$, $k\in\NNNzi$.
Denote by $\jet{k}{f}{x}$ the $k$-jet of $f$ at $x$.
Let $J^{k}(M,N)$ be the manifold of $k$-jets from $M$ to $N$ and
$d$ be a metric on $J^{k}(M,N)$.
For any $k\in\NNNzi$, the space $\smsp{M,N}$ can be
endowed with the ``weak'' and ``strong'' Whitney topologies,
which we denote by $C^{k}_{W}$ and $C^{k}_{S}$, resp.
(see e.g. Mather~\cite{Mather:IS2}, Hirsch~\cite{Hirsch:DiffTop}.)
Let us recall the definitions.

Let $\amap\in\smsp{M,N}$.
Then a base of the {\em weak\/} $C^{k}_{W}$-topology on $\smsp{M,N}$ at $\amap$
consists of sets of the form
\begin{equation}\label{equ:weak_base_nbh}
\NbhW{\amap}{K}{\eps}=
    \{ \bmap \in \smsp{M,N} \ | \
    d(\jet{k}{\amap}{x},\jet{k}{\bmap}{x}) < \eps, \forall x\in K\},
\end{equation}
where $K\subset\manif$ is compact and $\eps>0$.
For $k=0$, this topology is often called {\em compact-open}.
A base of the {\em strong\/} $C^{k}_{S}$-topology on $\smsp{M,N}$ at $\amap$
is generated by sets of the form
\begin{equation}\label{equ:strong_base_nbh_with_metric}
\NbhS{\amap}{\delta}=
 \{ g \in C^0(A,B) \ | \ d\left(\jet{k}{g}{x}, \jet{k}{f}{x} \right) < \delta(x),
            \forall x \in M \},
\end{equation}
where $\delta: M \to (0, \infty)$ is any continuous function.

Suppose $A$, $B$, $A'$, and $B'$ are manifolds,
$\XX\subset\smsp{A,B}$, $F:\XX \to \smsp{A',B'}$ is a map,
$r, r'\in\NNNzi$, and the symbols $T$ and $T'$
stand either for ``$W$'' or ``$S$''.
We say that $\XX$ is $C^{r}_{T}$-{\em open} ({\em -closed}, etc) if
it is open (closed) in the $C^{r}_{T}$-topology of $\smsp{A,B}$.
We say that $F$ is {\em $C^{\,r,r'}_{\,T,T'}$-continuous}
({\em -homeomorphism}, {\em -embedding}) if $F$ becomes
continuous (a homeomorphism, an embedding)
whenever $\smsp{A,B}$ and $\smsp{A',B'}$ are endowed with
the topologies $C^{\,r}_{\,T}$ and $C^{\,r'}_{\,T'}$, respectively.
If $r=r'$ and $T=T'$, then $F$ is said to be {\em $C^{r}_{T}$-continuous}.
Typical examples of $C^{\,r}_{\,T}$-continuous maps
are given, e.g., in Mather~\cite{Mather:IS2}.

\subsection{Flows}
Let $\Odzflow\subset\manif$ be an open, connected set and
$\Interval =(-a,a) \subset\RRR$, where $a>0$.
A {\em partial flow\/} on $\Odzflow$ is
a smooth map $\flow:\Odzflow\times\Interval \to \manif$
with the following properties.
If $x\in\Odzflow$ and $t,s\in\Interval$, then
\begin{enumerate}
\item
 $\flow(x,0) = x$,
\item
 $\flow( \flow(x,t), s) = \flow(x,t+s)$ provided $\flow(x,t)\in\Odzflow$
and $t+s\in\Interval$.
\end{enumerate}
If $\Odzflow=\manif$ and $\Interval=\RRR$,
then we say that $\flow$ is {\em global}.
By $\flow_t$ denote the restriction of $\flow$ to $\Odzflow \times t$,
where $t\in\Interval$.
The {\em trajectory\/} of a point $x\in\Odzflow$ is the set
$\flow(x\times\Interval)$.
There are three types of trajectories:
{\em constant\/} (fixed point),
{\em closed \em or \em periodic\/} (homeomorphic image of $S^1$),
and {\em non-closed} (one-to-one image of $\RRR$).
A point that is not fixed is called {\em regular}.
The {\em period\/} of a periodic point $x$ is the least positive number
$\Per\,x$ such that $\flow(x,\Per\,x)=x$.

The {\em Jordan cell} $J_p(A)$ of a $(k \times k)$-matrix $A$
is the following $(p k \times p k)$-matrix:
$$
J_p(A) =
\left.
\left\|
\begin{matrix}
A     & 0      & \cdots & 0 & 0 \\
E_k   & A      & \cdots & 0 & 0\\
\cdot & \cdot  & \cdot  & \cdot & \cdot \\
0 & 0 & \cdots & E_k    & A
\end{matrix}
\right\|
\ \right\} p,
$$
where $E_k$ is the unit $(k \times k)$-matrix.
If $\alpha, \beta\in\RRR$, then we denote
\begin{equation}\label{equ:rotation_matrix}
\Rot{\alpha,\beta} = \matr{\alpha}{-\beta}{\beta}{\alpha}.
\end{equation}
It is well known that each square matrix with real entries
is a conjugate to a matrix of the form
\begin{equation}\label{equ:Real_JordForm}
\diag[ J_{k_1}(\lambda_1), \ldots, J_{k_m}(\lambda_m),
J_{p_1}(\Rot{\alpha_1, \beta_1}), \ldots, J_{p_s}(\Rot{\alpha_s, \beta_s}) ],
\end{equation}
where $\lambda_i \in \RRR$, $(i=1...m)$ and $\beta_j\not=0$ $(j=1...s)$
(see e.g. Theorem~2.2.5 in Palis~J. and de Melo~W.~\cite{PalisMelo}).
We also need the following formulas
\begin{equation}\label{equ:exp_R}
e^{\Rot{\alpha, \beta} t} = e^{\alpha t}
\matr{\cos (\beta t)}{- \sin (\beta t)}{\sin (\beta t)}{\cos (\beta t)},
\end{equation}
\begin{equation}\label{equ:exp_J_A}
e^{J_{k}(A) t} =
\left\|
\begin{array}{cccc}
e^{At}                               & 0                                    & \cdots  & 0      \\
t \cdot e^{At}                       & e^{At}                               & \cdots  & 0      \\
\cdots                               & \cdots                               & \cdots  & \cdots \\
\frac{t^{k-1}}{(k-1)!}  \cdot e^{At} & \frac{t^{k-2}}{(k-2)!}  \cdot e^{At} & \cdots  & e^{At}
\end{array}
\right\|.
\end{equation}

\section{Maps generated by smooth shifts.}\label{sect:map_gen_flow}
Let $\flow:\Odzflow\times\Interval \to \manif$ be a partial flow.
Then $\flow$ yields the {\em shift\/} mapping $\Shift: \sectflow \to \shiftflow$
defined by
\begin{equation} \label{equ:h_flow}
\Shift(\afunc)(z) = \flow(z, \afunc(z)),
\end{equation}
where $\afunc \in \sectflow$ and $z \in \manif$.
If $\afunc\in\sectflow$,
then we say that $\amap=\Shift(\afunc)$ is a {\em shift along trajectories of $\flow$
by $\afunc$\/} and $\afunc$ is a {\em shift-function\/} of $\amap$.
\begin{lem}\label{lm:WWSSinvsh}
For each $\smcl\in\NNNzi$ and $T=$``$W$'' or ``$S$'',
the map $\Shift$ is $C^{\smcl}_{T}$-continuous.
\end{lem}
\proof
Let $*:\sectflow\to\inb\in\shiftflow$ be the constant map,
where $\inb:\Odzflow \subset \manif$ is the identity embedding.
Then $\Shift$ coincides with the following composition:
$$
\begin{CD}
\sectflow \xrightarrow{*\times\id}
\shiftflow \times  \sectflow \xrightarrow{\approx}
\smsp{\Odzflow, \manif \times \Interval} \xrightarrow{\flow_{*}}
\shiftflow,
\end{CD}
$$
where the first arrow is the product of $*$ and the identity map of $\sectflow$,
the second one is a natural $C^r_T$-homeomorphism,
and the third one is induced by $\flow$ and is
$C^r_T$-continuous as well (e.g. Mather~\cite{Mather:IS2}).
\qed

\begin{prop}\label{prop:gr_semigr}
If $\flow$ is global,
then the image $\IM\Shift$ is a semigroup and the intersection
$\IM\Shift \cap \DiffM$ is a subgroup of $\smmap$.
\end{prop}
\proof
Suppose $\afunc, \bfunc, \cfunc \in \smfunc$,
$\amap=\Shift(\afunc)$, $\bmap=\Shift(\bfunc)$,
$\cmap=\Shift(\cfunc)$, and $\cmap \in \DiffM$.
Let us show that $\amap \circ \bmap$ and $\cmap^{-1} \in \IM\Shift$.
Consider the functions
\begin{gather}
\label{equ:sh_gh}
\stfunc_{\amap \circ \bmap} = \bfunc + \afunc \circ \bmap, \\
\label{equ:sh_g_1}
\stfunc_{\cmap^{-1}} = - \cfunc \circ \cmap^{-1}.
\end{gather}
Then it can easily be seen that
$\Shift\left( \stfunc_{\amap \circ \bmap} \right) = \amap \circ \bmap$
and $\Shift\left( \stfunc_{\cmap^{-1}} \right) = \cmap^{-1}$.
\qed

\begin{defn}
The set $\Zid(\Shift) = \Shift^{-1}(\id_{\Odzflow})$,
where $\id_{\Odzflow}:U \subset \manif$ is the natural inclusion,
will be called the {\em kernel\/} of $\Shift$.
\end{defn}
Notice that, by formulas~\eqref{equ:sh_gh} and~\eqref{equ:sh_g_1},
$\Shift$ is not a homomorphism.
Nevertheless the following lemma explains our terminology.

\begin{lem}\label{lm:Zid_subgr}
Let $\afunc,\bfunc\in\sectflow$ be functions such that $\afunc-\bfunc \in \sectflow$.
Then $\Shift(\afunc)=\Shift(\bfunc)$ if and only if $\afunc-\bfunc \in \Zid$.
\end{lem}
\proof
The relation $\Shift(\afunc)(x)=\flow(x, \afunc(x)) = \flow(x, \bfunc(x))=\Shift(\bfunc)(x)$
for $x\in\Odzflow$ is equivalent to the following one:
$\flow(x, \afunc(x) - \bfunc(x)) = x$, i.e., $\afunc - \bfunc \in \Zid$.
\qed

\begin{cor}\label{cor:descript_Zid}
Suppose $\flow$ is global.
Then $\Zid=\Shift^{-1}(\id_{\manif})$ is a {\em subgroup\/} of $\smfunc$
and $\Shift^{-1}\circ \Shift(\afunc) = \{\afunc + \Zid\}$
for each $\afunc\in\smfunc$.
Thus there is a natural identification $\IM\Shift \cong \smfunc/\Zid$. \qed
\end{cor}

\subsection{Regular points}
Let $\pntb\in\Odzflow$ be a regular point of $\flow$,
$\afunc\in\sectflow$, $\amap=\Shift(\afunc)$, and $a=\afunc(\pntb)$.
Then the point $\pntc = \amap(\pntb)=\flow_{a}(\pntb)$ is also regular.
Hence there exists a neighborhood $\nbc$ of $\pntc$
and local coordinates $(x_1,\ldots,x_k)$ on $\nbc$
such that $\pntc=0$ and
$\flow(x_1,\ldots,x_k,t) = (x_1 + t, x_2,\ldots,x_k)$.
Let $\nbb = \amap^{-1}(\nbc) \cap \flow_{-a}(\nbc)$
be a neighborhood of $\pntb$.
Then $\afunc$ can be expressed on $\nbb$ in terms of $\amap$ as follows:
\begin{equation}\label{equ:reg_pnt_lsect}
  \afunc(x) = p_1\circ\amap\circ\flow(x,a) - p_1\circ\flow(x,a),
\end{equation}
where $p_1:\RRR^n \to \RRR$ is the projection onto the first coordinate.
This relation will be often used.
The first application is the local uniqueness of functions $\mu \in \Zid$ at regular points of $\flow$.
(Corollary~\ref{cor:reg_per_comp}.)
\begin{lem}\label{lm:mu_const_traj}
Let $\mu\in\Zid$, and let $\omega$ be a non-constant trajectory of $\flow$.
If $\omega$ is non-closed, then $\mu|_{\omega}=0$.
Otherwise, $\mu|_{\omega}= n\theta$, where $\theta=\Per\omega$ and $n \in \ZZZ$.
\end{lem}
\proof
If $\omega$ is non-closed,
then for any $x,y\in\omega$ there exists a {\em unique} number $t \in \Interval$
such that $y =\flow(x,t)$.
In particular, $t=0$ iff $x=y$. Thus $\mu|_{\omega}\equiv 0$.

Let $\omega$ be closed, and $x\in\omega$.
Then the relation $x =\flow(x,t)$ holds iff $t = n\theta$ for some $n\in\ZZZ$.
Since $\mu$ is continuous and the set $\{n\theta\}$ is discrete, it follows that
$\mu|_{\omega}$ is constant.
\qed

\begin{cor}\label{cor:reg_per_comp}
Let $\rcomp$ be a component of the set of regular points of $\flow$, $x\in\rcomp$ be a point,
and $\afunc,\bfunc\in\sectflow$ be such that $\Shift(\afunc)=\Shift(\bfunc)$.
If $\afunc(x)=\bfunc(x)$, then $\afunc|_{\rcomp} = \bfunc|_{\rcomp}$.
In particular, if $\afunc\in\Zid$ and $\afunc(x)=0$, then $\afunc|_{\rcomp}=0$.
\end{cor}
\proof
See formula~\eqref{equ:reg_pnt_lsect}.
\qed

\subsection{Fixed points}
\begin{prop}\label{pr:x_in_IntFix}
Let $\afunc, \bfunc \in\sectflow$ be such that
$\afunc = \bfunc$ on $\Odzflow\setminus\Int\Fix\flow$.
Then $\Shift(\afunc)=\Shift(\bfunc)$.
In particular, $\afunc\in\Zid$ iff $\bfunc\in\Zid$.
\end{prop}
\proof
We must show that
$\Shift(\afunc)(y)=\flow(y,\afunc(y))=\flow(y,\bfunc(y))=\Shift(\bfunc(y))$
for all $y\in\Odzflow$.
For $y\in\Odzflow\setminus\Fix\flow$
this holds by the condition $\afunc(y)=\bfunc(y)$.
If $y\in\Fix\flow$, then $\flow(y,t)=y$ for all $t\in\RRR$.
Hence $\flow(y,\afunc(y))=\flow(y,\bfunc(y))=y$.
\qed


\begin{prop}\label{pr:x_in_FrFix}
Let $\pnt\in\Fr(\Fix\flow)$
and $\nbhcomp = \cup_{\nind\in\Lambda} \nbhcomp_\nind$ be
the union of components $\nbhcomp_\nind$ of $\Odzflow\setminus\Fix\flow$
such that $\pnt \in \overline{\nbhcomp_\nind}$ for all $\nind\in\Lambda$.
If $\mu\in\Zid$,
then each of the following conditions \emph{(1)} and \emph{(2)}
implies that $\mu \equiv 0$ on $\overline{\nbhcomp}$.
\begin{enumerate}
\item[\emph{(1)}]
The tangent linear flow at $\pnt$ is trivial, i.e.,
$\frac{\partial\flow}{\partial x}(\pnt,t)=E_n$ for all $t\in\Interval$.
In particular, this holds for any point $\pnt \in \Fr(\Int\Fix\flow)$.
\item[\emph{(2)}]
$\mu(\pnt)=0$.
\end{enumerate}
\end{prop}
For the proof we need the following lemma.
Let $M(n)$ be the space of real square $n \times n$ matrices,
$E_n$ be the unit $n \times n$ matrix,
and $\exp:M(n)\to\GL_{n}(\RRR)$ be the exponential map.

\begin{lem}\label{lm:linear_map_periods}
Let $\{ A_i \}_{i\in\NNN} \subset M(n)$ be a sequence of matrices
such that for each $i\in\NNN$ the flow $\flow_i(x,t) = e^{A_i t}x$
has at least one closed trajectory.
Let $\theta_i$ be the minimum of periods of the closed trajectories of $\flow_i$.
If $\lim\limits_{i \to \infty} A_i = 0$, then
$\lim\limits_{i \to \infty} \theta_i = \infty$.
\end{lem}
\proofstyle{Sketch of proof.}
Let $A \in M(n)$ and $\Lambda = \{ \lambda_k \}_{k=1}^{r}$
be the set of eigenvalues of $A$ such that
$\lambda_k = i \beta_k$ for $\beta_k \in \RRR\setminus\{0\}$.
Then it is easily seen that the linear flow $\flow(x,t)=e^{At}x$
has a closed trajectory iff $\Lambda \not= \emptyset$.
In this case the period of any closed trajectory of $\flow$ is \
$\geq \min\limits_{k=1...r} \frac{2\pi}{|\beta_k|} $.
Now the lemma follows from the continuity of spectrums
of linear operators.
\qed

\proofstyle{Proof of Proposition~\ref{pr:x_in_FrFix}.}
We will show that $\mu \equiv 0$ on each component $\nbhcomp_\nind$.
Suppose that $\nbhcomp_\nind$ contains a non-closed trajectory.
Then $\mu|_{\nbhcomp_\nind} \equiv 0$ by Corollary~\ref{cor:reg_per_comp}.
Thus suppose that $\nbhcomp_\nind$ consists of periodic points only.
Let $\{ \pnt_i\}_{i\in\NNN} \subset \nbhcomp_\nind$
be a sequence of periodic points of $\flow$ converging to $\pnt$ as $i\to\infty$.
For each $i\in\NNN$ denote $\theta_i=\Per\,\pnt_i$.
By Lemma~\ref{lm:mu_const_traj}, $\mu(\pnt_i) = n_i\theta_i$ for some $n_i\in\ZZZ$.
Then, by continuity of $\mu$, we get $\mu(\pnt_i) = n_i \theta_i \to \mu(\pnt) < \infty$.
Taking a subsequence of $\{ \theta_i\}_{i\in\NNN}$ (if necessary) we can assume
that there exists a finite or an infinite limit $\theta = \lim\limits_{i \to \infty} \theta_i \geq 0$.

We prove that in both cases (1) and (2),
$n_i=0$ and therefore $\mu(\pnt_i)=0$ for all sufficiently large $i\in\NNN$.
Since each point $\pnt_i\in\nbhcomp_\nind$ is regular,
it follows that $\mu\equiv0$ on $\nbhcomp_\nind$.
This will complete our proposition.

Define $\dflow: \Odzflow \times \Interval \to \GL_{n}(\RRR)$
by $\dflow(x,t) = \frac{\partial \flow}{\partial x}(x,t)$.
Since $\dflow(x,0)=E_n$ for all $x \in \Odzflow$,
it follows that the map
$\expflow = \exp^{-1} \circ \dflow:\Odzflow \times \Interval \to M(n)$
is well defined on some neighborhood of $(\pnt,0)$ in $\Odzflow \times \Interval$.
Thus $\dflow(x,t) = e^{\expflow(x,t)}$.

Notice that the restriction map $\dflow(x,*):\Interval \to \GL_{n}(\RRR)$
is a local homomorphism, i.e., it yields a linear flow on $\RRR^n$.
Hence $A(x,t) = A(x) = \expflow(x,t)/t$ does not depend on $t\in\Interval$
and $\dflow(x,t) = e^{A(x) t}$.

We now show that for each periodic point $x$ the flow $\dflow(x,*)$ has closed trajectories.
Let $\Fld(x) = \frac{\partial\flow}{\partial t}(x,0)$ be the vector field generating $\flow$.
Applying the operator $\frac{\partial}{\partial t}$ to both parts of the relation
$\flow( \flow(x,t), s) = \flow( \flow(x,s), t)$ and then setting $s=0$, we obtain
$$
\begin{CD}
\frac{\partial\flow}{\partial t} (\flow(x,t), 0) =
\frac{\partial\flow}{\partial x} (x, t) \frac{\partial\flow}{\partial t} (x, 0),
\end{CD}
$$
i.e., $\Fld(\flow(x,t)) = \dflow(x,t) \Fld(x)$. Thus the vectors
$\Fld(\flow(x,t))$ and $\Fld(x)$ lie on same trajectory of the flow $\dflow(x,*)$.
In particular, if $x$ is a periodic point of $\flow$,
then $\dflow(x, \Per\,x)\Fld(x) = \Fld(x)$, i.e.,
the vector $\Fld(x)$ is a periodic point of $\dflow(x,*)$ and $\Per\Fld(x) \leq \Per(x)$.
Consider now the cases (1) and (2) of our proposition.

(1)
Suppose $\frac{\partial\flow}{\partial x}(\pnt,t) = \dflow(\pnt,t) = E_n$.
Since $\dflow(\pnt_i,t) \to \dflow(\pnt,t) = E_n$,
it follows from Lemma~\ref{lm:linear_map_periods} that
$\theta=\lim\limits_{i \to \infty} \theta_i \geq \lim\limits_{i \to \infty} \Per\,\Fld(\pnt_i) = \infty$.
Since $\lim\limits_{i \to \infty} n_i\theta_i < \infty$, we get
$\lim\limits_{i \to \infty} n_i = 0$.

(2)
Let $\mu(\pnt)=0$.
We claim that $\theta = \lim\limits_{i \to \infty} \theta_i > 0$.
Since $\lim\limits_{i \to \infty} n_i\theta_i = \mu(z) = 0$,
it will follow that $\lim\limits_{i \to \infty} n_i = 0$.
Thus suppose $\theta=0$. Then
$$
\begin{CD}
\frac{\partial\flow}{\partial x}(\pnt, t) @<<{i\to\infty}<
\frac{\partial\flow}{\partial x}(\pnt_i, t) =
\frac{\partial\flow}{\partial x}(\pnt_i, \theta_i \{ t / \theta_i \} ) @>>{i\to\infty}>
\frac{\partial\flow}{\partial x}(\pnt_i, 0) = E_n,
\end{CD}
$$
i.e., $\frac{\partial\flow}{\partial x}(\pnt, t) = E_n$ for all $t \in \Interval$,
where $\{t\}$ is the fractional part of $t \in \RRR$.
Then by (1) we get $\theta = \infty$, which contradicts the assumption $\theta=0$.
\qed

\begin{thm}\label{th:Zid_descr}
\begin{enumerate}
\item [\emph{(1)}]
Suppose $\Int\Fix\flow \not= \emptyset$. Then
$$\Zid = \{ \mu \in \sectflow \ | \ \mu|_{\Odzflow \setminus \Int\Fix\flow} = 0\}.$$
\item [\emph{(2)}]
Let \,$\Int\Fix\flow = \emptyset$.
Then either $\Zid=\{0\}$ or there exists a function $\nu\in\Zid$
such that $\nu(x)\not=0$ for all $x\in\Odzflow$ and
for any other $\mu\in\Zid$, we have $\mu = n \nu$, where $n\in\ZZZ$.
Thus if $\flow$ is global, then $\Zid$ is either $0$ or $\ZZZ$.
\end{enumerate}
\end{thm}
\proof
(1)
Let $X = \Odzflow \setminus \Int\Fix\flow$ and $\mu\in\sectflow$.
If $\mu|_{X}=0$, then by Proposition~\ref{pr:x_in_IntFix}, we have $\mu \in \Zid$.

Conversely, let $\mu\in\Zid$.
We will show that the set $Y = \mu^{-1}(0) \cap X$ is a nonempty,
open-closed subset of $X$ intersecting each component of $X$.
This implies $\mu|_{X}=0$ whence $Y=X$.

Evidently, $Y$ is closed.
Moreover, by (1) of Proposition~\ref{pr:x_in_FrFix}, we have that  $\mu|_{\Fr(\Int\Fix\flow)} = 0$.
Thus $Y \supset \Fr(\Int\Fix\flow) = \Fr(X) \not= \emptyset$.
Since $\Odzflow$ is connected,
it follows that $\Fr(\Int\Fix\flow)$ (and therefore $Y$) intersects each component of $X$.
Finally, let $\pnt\in Y$.
Then, from Corollary~\ref{cor:reg_per_comp}, we obtain that
$\mu=0$ in some neighborhood of $\pnt$ in $X$ whenever $\pnt$ is regular.
Suppose $\pnt\in\Fix\flow$.
Then by (2) of Proposition~\ref{pr:x_in_FrFix},
$\mu = 0$ in some neighborhood of $\pnt$ in $X$.
Thus $Y$ is open in $X$.

(2)
Let $\Int\Fix\flow = \emptyset$.
For each $x \in \Odzflow$ define $\tau_x:\Zid \to \RRR$ by $\tau_x(\mu)=\mu(x)$, where $\mu\in\Zid$.
Then by Corollary~\ref{cor:reg_per_comp} and
by (2) of Proposition~\ref{pr:x_in_FrFix},  $\tau_x$ is injective.
Moreover, it follows from Lemma~\ref{lm:mu_const_traj} that for each regular point $x$ of $\flow$
the set $\IM\tau_x$ is a discrete subset of $\RRR$
consisting of some integer multiples of some $\theta\in\RRR$.
Suppose $\Zid\not=\{0\}$. Since $\Interval$ is a connected neighborhood of $0\in\RRR$,
we see that there exists a number $g \in \IM\tau_x$ dividing all elements of $\IM\tau_x$.
Then the function $\nu = \tau_x^{-1}(g)$ satisfies the statement of the theorem.
\qed

\begin{prop}\label{pr:suffcond_Zid_zero}
Suppose $\Int\Fix\flow = \emptyset$.
Then each of the following two conditions implies that $\Zid=\{0\}$.
\begin{enumerate}
\item[\emph{(1)}]
The tangent flow at some fixed point $\pnt\in\Fix\flow$ is trivial;
\item[\emph{(2)}]
$\flow$ has at least one non-closed trajectory.
\end{enumerate}
\end{prop}
\proof
Let $\mu\in\Zid$.
From (1) of Proposition~\ref{pr:x_in_FrFix} and
Lemma~\ref{lm:mu_const_traj} it follows that
both conditions (1) and (2)
imply that there exists a point $x\in\Odzflow$ such that $\mu(x)=0$.
Since $\Int\Fix\flow = \emptyset$, it follows from Theorem~\ref{th:Zid_descr}
that $\mu \equiv 0$ on $\Odzflow$.
\qed


\section{Local sections of $\Shift$}\label{sect:loc_sect_of_invsh}
Suppose that $\flow$ is global.
Then by Corollary~\ref{cor:descript_Zid}, the map $\Shift$ has the following decomposition:
\begin{equation}\label{equ:decomp_invsh}
\Shift : \smfunc \xrightarrow{\widetilde{\Shift}} \smfunc/\Zid
\xrightarrow{j} \IM\Shift \subset \smmap,
\end{equation}
where $\widetilde\Shift$ is the factor-map and $j$ is the bijection $\{\afunc+\Zid\}\mapsto\Shift(\afunc)$.

Suppose $\smfunc$ and $\smmap$ are endowed with some topologies.
Recall that the corresponding factor-topology on
$\smfunc / \Zid$ is defined as follows:
a subset $\nbhtmp \subset \smfunc/\Zid$ is open iff
$\widetilde\Shift^{-1}(\nbhtmp)$ is open in $\smfunc$.
Then $j$ is continuous iff so is $\Shift$.
In general, $j$ is not a homeomorphism, i.e., the factor-topology of $\IM\Shift$
from $\smfunc$ can differ from the induced topology of the ambient space $\smmap$.
For the case $\Int\Fix\flow=\emptyset$
we give a sufficient condition for $j$ to be a homeomorphism
in the related strong Whitney topologies.
It requires existence of weakly continuous local sections of $\Shift$
at each fixed point of $\flow$.
For purposes of \MainTheorem, we will consider flows depending on a parameter.

Let $\flow$ be a partial flow and $\Diskk$ be an open $k$-dimensional disk.
Define the partial flow
$\tflow:(\Odzflow\times\Diskk)\times\Interval\to\manif\times\Diskk$
as the product of $\flow$ by the trivial flow on $\Diskk$, i.e.,
$\tflow(x,\tau,t)=(\flow(x,t), \tau)$,
where $(x,\tau,t)\in\Odzflow\times\Diskk\times\Interval$.
For each open subset $\nbs \subset \Odzflow \times \Diskk$
formula~\eqref{equ:h_flow} yields a {\em shift-map\/}
$\ShM{\nbs}:\smsp{\nbs,\Interval} \to \smsp{\nbs,\manif}$
of $\tflow$ by $\ShM{\nbs}(\afunc)(x,\tau) = \flow(x,\tau, \afunc(x,\tau))$,
where $(x,\tau)\in V$ and $\afunc \in \smsp{\nbs,\manif}$.
\begin{prop}\label{prop:inject_phiV}
The map $\Shift$ is {\em locally injective\/} in $C^{\smcl}_{S}$-topology of $\sectflow$
for any $\smcl\in\NNNzi$ iff \/ $\Int\Fix\flow=\emptyset$.
In this case there exists a continuous function $\unifunc:\Odzflow\to(0,\infty)$
such that for any open $\nbs\subset\Odzflow\times\Diskk$ and $\afunc\in\smsp{\nbs,\Interval}$
the restriction of $\ShM{\nbs}$ to the $C^{0}_{S}$-neighborhood
\begin{equation}\label{equ:Nbf_unifunc}
 \NBfdnbs = \{ \bfunc\in\smsp{\nbs,\Interval} \ | \
      |\afunc(x,\tau)-\bfunc(x,\tau)|< \unifunc(x), \forall (x,\tau)\in\nbs \}
\end{equation}
of $\afunc$ is injective and
$\NBfdnbs \cap \{ \NBfdnbs + \mu \} = \emptyset$ for each $\mu\in\Zid$ provided $\mu\not=0$.
Hence, if $\flow$ is global, then $\widetilde{\Shift}$ is covering map.
\end{prop}
\proof
Suppose $\Int\Fix\flow\not=\emptyset$.
Let $\afunc \in \sectflow$ and let $\NN$ be any $C^{\smcl}_{S}$-neighbor\-hood of $\afunc$
in $\sectflow$.
Then there exists $\bfunc\in\NN$ such that $\afunc=\bfunc$ on
$\manif\setminus\Int\Fix\flow$, whence $\Shift(\afunc)=\Shift(\bfunc)$
by Lemma~\ref{pr:x_in_IntFix}. Thus $\Shift$ is not injective.

Let $\Int\Fix\flow=\emptyset$.
We will construct a function $\unifunc$ satisfying the statement of our proposition.
From the proof of Proposition~\ref{pr:x_in_FrFix} it follows that
for each $x\in\Odzflow$ there exists a compact neighborhood $\nbhtmp_x$ of $x$
and a number $\tau_x>0$ such that $\tau_x$
is less than the half of the period of any closed trajectory of $\flow$
passing through $\nbhtmp_x$.
If $\nbhtmp_x$ intersects no periodic trajectories, then we set $\tau_x=1$.
Since $\manif$ is paracompact,
there exists a continuous function $\unifunc:\Odzflow\to(0,1)$
such that $\unifunc(x)<\tau_x$ for all $x\in\Odzflow$.
Then $\delta$ satisfies the statement of our proposition.
\qed

\begin{defn}\label{def:PLS}
Suppose $\Int\Fix\flow=\emptyset$.
A point $\pnt\in\Odzflow$ is said to be an {\em\condSk}-point of $\flow$
if for any sufficiently small
open neighborhood $\nbs\subset\Odzflow\times\Diskk$ of $(\pnt,0)$
with compact closure $\overline{\nbs}$,
and any function $\afunc\in\smsp{\nbs,\Interval}$
there exists a $C^{0}_{W}$-neighborhood $\NBf\subset\smsp{\nbs,\Interval}$
of $\afunc$
such that the restriction of $\ShM{\nbs}$ to $\NBf$ is injective
and the inverse map
$(\ShM{\nbs})^{-1}: \ShM{\nbs}(\NBf)\to\NBf$
is $C^{\smk}_{W}$-continuous for each $\smk\in\NNNz$.
A point is {\em\condS} if it is {\em\condSk} for each $k\in\NNNz$.
\end{defn}
Let us explain this definition in more details.
From Proposition~\ref{prop:inject_phiV}, we see that
for any open neighborhood $\nbs\subset\Odzflow\times\Diskk$ of $(\pnt,0)$
there exists a $C^{0}_{W}$-neighborhood $\NBf\subset\smsp{\nbs,\Interval}$
of $\afunc$ such that the restriction of $\ShM{\nbs}$ to $\NBf$ is injective
(e.g., we may put $\NBf=\NBfdnbs$).
Thus, on $\NBf$, the relation $\bmap(x,\tau)=\flow(x,\tau, \bfunc(x,\tau))$
is equivalent to $\bfunc(x,\tau) = \Psi(x,\tau,\bmap)$,
where $\Psi$ is some function.
Then a point $\pnt$ is \condSk\ whenever $\Psi$
induces a continuous map $\ShM{\nbs}(\NBf)\to\NBf$
in the {\em weak\/} $C^{\smk}$-topologies.
In particular, this holds whenever $\Psi$ is smooth in $(x,\tau)$
and continuously depend on all partial derivatives of $\bmap$ near $(\pnt,0)$
up to order $\smk$.
For instance, the following lemma is a direct corollary of formula~\eqref{equ:reg_pnt_lsect}.
\begin{lem}\label{lm:lsect_reg_pnt}
Any regular point of a flow is {\em\condS}.\qed
\end{lem}

\begin{thm}\label{th:LocSect}
Suppose $\Int\Fix\flow = \emptyset$ and
each fixed point $\pnt\in\Fix\flow$ of $\flow$ is {\em\condSz}.
For any $\afunc\in\sectflow$ let
$$
 \NBfd = \NBfd(\afunc) =
         \{ \bfunc\in\sectflow \ | \ |\bfunc(x)-\afunc(x)| < \unifunc(x), \forall x\in\Odzflow \}
$$
be a $C^{0}_{S}$-neighborhood of $\afunc$ in $\sectflow$ such that
the restriction $\Shift|_{\NBfd}$ is injective.
Then the inverse map $\Shift^{-1}: \Shift(\NBfd) \to \NBfd$
is $C^{\smk}_{S}$-continuous for any $\smk\in\NNNzi$.
Hence, for global $\flow$, the map $\Shift:\smfunc\to\IM\Shift$ is a covering map in
the $C^{\smk}_{S}$-topologies for any $\smk\in\NNNzi$.
\end{thm}
For the proof we need the following statement.
\begin{lem}\label{cor:intersect_of_preimages}
Let $M$ and $N$ be smooth manifolds, $\amap \in \smsp{M,N}$, and $\smk\in\NNNz$.
Let also $\{\nbi\}_{i\in\Lambda}$ be a locally finite family of open subsets of $M$
and for each $i\in\Lambda$ let $\UUi$ be a $C^{\smk}_{W}$-neighborhood of the restriction
$\amap|_{\nbi}$ in the space $\smsp{\nbi,N}$.
Define $p_i:\smsp{M,N} \to \smsp{\nbi,N}$ by $\amap \mapsto \amap|_{\nbi}$ for $\amap\in\smsp{M,N}$,
and let $\tUUi = p_i^{-1}(\UUi)$.
Then $\tUU = \mathop\cap\limits_{i\in\Lambda} \tUUi$
is a $C^{\smk}_{S}$-neighborhood of $\amap$ in $\smsp{M,N}$.
\end{lem}
\proof
Since $p_i$ is $C^{\smk}_{W}$-continuous,
it follows that for each $i\in\Lambda$
the set $\tUUi$ contains a base $C^{\smk}_{W}$-neighborhood
$\NbhW{\amap|_{\nbi}}{\Ki}{\epsi}$ of $\amap|_{\nbi}$
defined by formula~\eqref{equ:weak_base_nbh},
where $\Ki \subset \nbi$ is compact and $\epsi>0$.
Note that $\{\Ki\}_{i\in\Lambda}$ is locally finite and $M$ is paracompact.
So there exists a continuous function
$\delta:M \to (0,\infty)$ such that $\delta(x)<\epsi$ for $x \in \Ki$.
Let $\tNN = \NbhS{\amap}{\delta}$
be the open $C^{\smk}_{S}$-base neighborhood of $\amap$ in $C^{\smk}(M,N)$
defined by formula~\eqref{equ:strong_base_nbh_with_metric}.
Then, obviously, $\tNN \subset \tUU$.
\qed

\proofstyle{Proof of Theorem~\ref{th:LocSect}.}
From Lemma~\ref{lm:lsect_reg_pnt} and the assumption of the theorem we obtain
that each $\pnt\in\Odzflow$ is an \condSz-point of $\flow$.
Let $\smk\in\NNNz$, $d$ be a metric on $J^{\smk}(\Odzflow,\Interval)$,
and $\unifuncA:\Odzflow\to(0,\infty)$ be a continuous function
such that the base $C^{\smk}_{S}$-neighborhood
$
 \NBfdA= \{
               \bfunc\in\sectflow \ | \
                 d(\jet{\smk}{\afunc}{x}, \jet{\smk}{\bfunc}{x}) < \unifuncA(x),
                \forall x\in\Odzflow
            \}
$
of $\afunc$ (see formula~\eqref{equ:strong_base_nbh_with_metric})
lies in $\NBfd$.
Hence $\Shift|_{\NBfdA}$ is injective.
Let us prove that $\Shift(\NBfdA)$ contains
a $C^{\smk}_{S}$-neighborhood of $\amap=\Shift(\afunc)$
in $\IM\Shift$, i.e.,
there exists an open neighborhood $\tNBm$ of $\amap$
in $\shiftflow$ such that $\IM\Shift \cap \tNBm \subset \Shift(\NBfdA)$.
Since $\unifuncA$ can be chosen arbitrary small,
we will get that $\Shift^{-1}$ is $C^{\smk}_{S}$-continuous on $\Shift(\NBfd)$.

Since $\Odzflow$ is paracompact,
there exist two at most countable locally finite coverings
$\{ \nbi \}_{i\in\Lambda}$ and $\{ \Ki \}_{i\in\Lambda}$ of $\Odzflow$
such that $\nbi$ is open and $\Ki \subset \nbi$ is compact.
Moreover, using the \condSz\ property of points,
we may assume that
for each $i\in\Lambda$ there exists a $C^0_{W}$-neighborhood $\NBfnbi$ of $\afunc|_{\nbi}$
in $\smsp{\nbi,\Interval}$ such that the restriction of $\ShM{\nbi}$ to $\NBfnbi$
is a $C^{\smcl}_{W}$-embedding for all $\smcl\in\NNNz$.

Let $d_i$ be a metric on $J^r(\Odzflow,\Interval)$ and $\epsi = \inf\limits_{x\in\Ki}\unifuncA(x)$.
Define a $C^{\smcl}_{W}$-neighborhood $\NbhW{\afunc|_{\nbi}}{\Ki}{\epsi}$
of $\afunc|_{\nbi}$ by formula~\eqref{equ:weak_base_nbh}, with $d=d_i$,
and set $\tNBfnbi = \NbhW{\afunc|_{\nbi}}{\Ki}{\epsi} \cap \NBfnbi$.
Then the image $\NBmi = \ShM{\nbi}(\tNBfnbi)$ is a
$C^{\smcl}_{W}$-neigh\-bor\-hood of $\amap|_{\nbi}$ in
$\smsp{\nbi,\manif}$ for all $\smcl\in\NNNz$.

Let $\rMmapi:\shiftflow \to \smsp{\nbi,\manif}$
be the ``restriction to $\nbi$'' mapping
and $\tNBmi = \rMmapi^{-1}(\NBmi)$.
Since $\{\nbi\}_{i\in\Lambda}$ is locally finite,
it follows from Lemma~\ref{cor:intersect_of_preimages} that the intersection
$\mathop\cap\limits_{i\in\Lambda} \tNBmi$ contains some
$C^{\smk}_{S}$-open neighborhood
$\tNBm$ of $\amap$.
One can easily verify that $\tNBm \cap \IM\Shift \subset \Shift(\NBfdA)$.
This proves the theorem.
For the convenience of the reader we
give the following commutative diagram illustrating our
constructions.
$$
\begin{CD}
\NBfdA         \ & \ \subset \ & \ \sectflow  @>>>
    \smsp{\nbi,\Interval} \ & \ \supset \ & \ \tNBfnbi = \NbhW{\afunc|_{\nbi}}{\Ki}{\epsi} \cap \NBfnbi \\
@V{\Shift}VV @V{\Shift}VV @V{\ShM{\nbi}}VV @V{\ShM{\nbi}}VV \\
\Shift(\NBfdA) \ & \ \subset \ & \ \shiftflow @>{p_i}>>
    \smsp{\nbi,\manif} \ & \ \supset \ & \ \ShM{\nbi}(\tNBfnbi) = \NBmi.
\end{CD}
$$


\section{Shifts that are diffeomorphisms}\label{sect:diff_gen_flow}
The following Theorem~\ref{th:h_diff} gives a precise description of the set
$$
\ImFn = \Shift^{-1}(\DiffM) =
  \Shift^{-1}(\IM\Shift \cap \DiffM) \subset \smfunc.
$$

\begin{thm}\label{th:h_diff}
Suppose that $\flow$ is a global flow on $\manif$,
$\Fld(\pnt) = \frac{\partial \flow}{\partial t}(\pnt, 0)$
be the vector field generating $\flow$ and $\afunc \in \smfunc$.
Then $\amap = \Shift(\afunc)$ is a diffeomorphism of $\manif$ if and only if
$\amap$ is proper ($\amap^{-1}(K)$ is compact for each compact $K\subset\manif$)
and the following inequality holds at each $z \in \manif$:

\begin{equation}\label{equ:ds/dF<>-1}
d\afunc(\Fld(\pnt)) \not= -1.
\end{equation}
Moreover, $\amap$ preserves orientation of trajectories iff \ $d\afunc(\Fld(\pnt)) > -1$.
\end{thm}
\proof
The necessity is implied by the following lemma.
\begin{lem}\label{lm:shift_loc_diff}
The tangent map $\amap'(\pnt):\tanpm\to\tanpm$ is an isomorphism
if and only if ~\eqref{equ:ds/dF<>-1} holds true at $z$.
\end{lem}
\newcommand\GLOBDET{\left|\begin{array}{cc} A_x + A_t\,\afunc_x & A_y + A_t\,\afunc_y \\ B_x + B_t\,\afunc_x & B_y + B_t\,\afunc_y \\ \end{array}\right|}
\newcommand\MAINDET{\left|\begin{array}{cc} A_x & A_y \\ B_x & B_y \\ \end{array}\right|}
\newcommand\DETX{\left|\begin{array}{cc} A_t & A_y \\ B_t & B_y \\ \end{array}\right|}
\newcommand\DETY{\left|\begin{array}{cc} A_x & A_t \\ B_x & B_t \\ \end{array}\right|}
\proof
We can assume $\afunc(\pnt)=0$.
Otherwise, set $\bfunc = \afunc - \afunc(\pnt)$ and $\bmap = \Shift(\bfunc)$.
Then $\bfunc(\pnt)=0$, $d\bfunc(\Fld)=d\afunc(\Fld)$,
and $\bmap'(\pnt)$ is an isomorphism iff so is $\amap'(\pnt)$.

So let $\afunc(\pnt)=0$. Then $\amap(\pnt)=\pnt$.
Let us choose a local chart at $\pnt$ such that $\pnt=0\in\RRR^n$
and calculate the determinant $|\amap'(\pnt)|$ of $\amap'(\pnt)$.
We claim that
\begin{equation}\label{equ:dh_0}
|\amap'(\pnt)| = 1 + d\afunc(\Fld(\pnt)).
\end{equation}
This will prove our lemma.
For simplicity we consider only the case $n=2$.
The general case is analogous.
Let $\flow = (A, B)$, where $A(x,y;t)$ and $B(x,y;t)$ are coordinate
functions of $\flow$ at this chart. Then
\begin{multline*}
|\amap'(\pnt)|  =  |\flow'(\pnt; \afunc(\pnt))|  =  \GLOBDET = \\
 = \MAINDET + \DETX\ \afunc_x + \DETY\ \afunc_y
 = \MAINDET \left( 1 + X \afunc_x + Y \afunc_y \right),
\end{multline*}
where, by Cramer's formulas, the vector $(X,Y)$ is a solution of the following linear equation
$$
\matr{A_x}{A_y}{B_x}{B_y}
\left\| \begin{array}{c} X \\ Y \end{array} \right\| =
\left\| \begin{array}{c} A_t \\ B_t \end{array} \right\|.
$$
Since $\afunc(\pnt)=0$, we have
$\flow_{\afunc(\pnt)} =\matr{A_x}{A_y}{B_x}{B_y}=\id$
and $(A_t, B_t) = \Fld(\pnt)$, whence $(X,Y)=\Fld(\pnt)$.
It remains to note that \
$X \afunc_x + Y \afunc_y = d\afunc(X,Y) = d\afunc(\Fld(\pnt))$.
\qed

\proofstyle{Sufficiency.}
Suppose $\amap$ is proper and~\eqref{equ:ds/dF<>-1} holds
at each $\pnt\in\manif$.
Then $\amap'$ is non-degenerate at each $\pnt\in\manif$
and it remains to show that
$\amap$ is bijective.
Since $\amap(\omega) \subset \omega$ for each trajectory $\omega$ of $\flow$,
we should establish that
$\amap(\omega)=\omega$ and $\amap|_{\omega}$ is one-to-one.
This is evident for fixed points.
Let $\omega$ be a regular trajectory.
Since $\amap'$ is non-degenerate, it follows that the restriction
$\amap|_{\omega}:\omega\to\omega$
is a proper map having no critical points and homotopic to a diffeomorphism.
Hence $\amap|_{\omega}$ is one-to-one.
\qed
\begin{cor}
\label{cor:zfix_hdiff}
Let $\pnt\in\Fix\flow$ and $\afunc\in\smfunc$.
Then the map $\Shift(\afunc)$ is a local diffeomorphism near $\pnt$.
\end{cor}
\proof
Since $\Fld(\pnt)=0$, we see that $d\afunc(\Fld(\pnt)) = 0>-1$.
\qed

\subsection{Decomposition of $\ImFn$}
Define the subsets $\ImFnplus$ and $\ImFnminus$ of $\smfunc$ by
\begin{gather*}
\ImFnplus = \{ \afunc \in \ImFn \ | \ d\afunc(F(x)) > -1 , \forall x\in\manif \}, \\
\ImFnminus = \{ \afunc \in \ImFn \ | \ d\afunc(F(x)) < -1, \forall x\in\manif \}.
\end{gather*}
Then $\ImFnplus \cap \ImFnminus = \emptyset$.
Moreover, from Theorem~\ref{th:h_diff} we get $\ImFn = \ImFnplus \cup \ImFnminus$.

\begin{lem}\label{lm:ImFnplus_minus_convex}
For each $\smcl\in\NNNzi$
the sets $\ImFnplus$ and $\ImFnminus$ are $C^{\smcl}_{S}$-open in $C^{\infty}(\manif)$
and {\em convex\/} if regarded as subsets of the {\em linear space\/} $C^{\infty}(\manif)$.
\end{lem}
\proof
Since $\DiffM$ is $C^{\smcl}_{S}$-open in $\smmap$,
and $\Shift$ and differentiating along the vector field
are $C^{\smcl}_{S}$-continuous, we obtain that $\ImFnplus$ and $\ImFnminus$ are also open.
Let us prove that $\ImFnplus$ is convex.
The proof for $\ImFnminus$ is analogous.

Let $\afunc_0, \afunc_1 \in \ImFnplus$, $\afunc_s=s \afunc_0 + (1-s) \afunc_1$,
and $\amap_{s} = \Shift(\afunc)$ for $s\in[0,1]$.
Then $d\afunc_s(F) = s\,d\afunc_0 + (1-s)\,d\afunc_1 > -1$, whence
$\amap'_s(x)$ is an isomorphism for each $x \in \manif$.
By the arguments similar to the proof of sufficiency in Theorem~\ref{th:h_diff},
$\amap_s$ is injective. Let us show that $\amap_s$ is onto.

First consider the flow
$\Shflow:(\manif\times\RRR)\times \RRR \to\manif\times\RRR$ on $\manif\times\RRR$
defined by the formula $\Shflow(x,t,s)=(x,t+s)$.
Then the {\em mapping\/} $\flow:\manif\times\RRR\to\manif$ gives rise to the factorization of
the flow $\Shflow$ onto the {\em flow\/} $\flow$, i.e.,
$\flow \circ \Shflow_s = \flow_s \circ \flow$ for all $s\in\RRR$. Indeed,
$$
  \flow \circ \Shflow_s(x,t) = \flow(x, t+s) = \flow(\flow(x,t),s)= \flow_s \circ \flow(x,t).
$$
Let $\afunc \in \smfunc$, $\amap(x)=\flow(x,\afunc(x))$,
$\tafunc=\afunc\circ\flow:\manif\times\RRR\to\RRR$, and
$\tamap:\manif\times\RRR\to\manif\times\RRR$ be the shift along trajectories of $\Shflow$
by the function $\tafunc$,
i.e., $\tamap(x,t)=\Shflow(x,t,\tafunc(x,t)) = (x, t + \afunc\circ\flow(x,t))$.
Then it is easily seen that
\begin{equation}\label{equ:CD_factor_traj}
\flow \circ \tamap = \amap \circ \flow.
\end{equation}
Now let us define $\tafunc_s=\afunc_s \circ \flow$ and $\tamap_s(x,t)=\Shflow(x,t, \tafunc_s(x,t))$.
Then
\begin{multline} \label{equ:tld_amap_t}
\tamap_s(x,t)=
(x,t + \tafunc_s(x,t)) =
(x,\,s(t + \tafunc_0(x,t))\, +\,(1-s)(t+\tafunc_1(x,t))\,).
\end{multline}
By assumption, $\amap_0$ and $\amap_1$ are onto.
Then, by~\eqref{equ:CD_factor_traj}, so are $\tamap_0$ and $\tamap_1$.
It follows from formula~\eqref{equ:tld_amap_t} that so is $\tamap_s$
for each $s\in[0,1]$. Hence, by~\eqref{equ:CD_factor_traj}, $\amap_t$ is onto.
\qed

Let us illustrate Lemmas~\ref{lm:shift_loc_diff} and~\ref{lm:ImFnplus_minus_convex}
by applying them to the flow $\flow(x,t) = x+t$ on $\RRR$.
Let $\afunc\in\smsp{\RRR,\RRR}$ and $\amap(x)=\flow(x,\afunc(x)) = x + \afunc(x)$.
Then the Lie derivative of $\afunc$ with respect to $\flow$ coincides with the usual derivative $\afunc'$.
Hence the inequality $\afunc'(z)\not=-1$ is equivalent to $\amap'(z)\not=0$
(compare with Lemma~\ref{lm:shift_loc_diff}).
Notice that $\RRR$ is a unique trajectory of $\flow$.
So the space of diffeomorphisms of $\RRR$ preserving
(reversing) the orientation of this trajectory
coincides with the space of smooth monotone functions
with positive (negative) derivative.
Therefore, these spaces are open and convex (compare with Lemma~\ref{lm:ImFnplus_minus_convex}).


\section{Describing the image $\IM\Shift$} \label{sect:image_of_invsh}
Let $\flow$ be a {\em global} flow on $\manif$.
Define the subset $\Cflow \subset \smmap$ to be consisting of maps $\amap$ such that
\begin{enumerate}
\item $\amap(\omega) \subset \omega$ for each trajectory $\omega$ of $\flow$,
\item $\amap'(\pnt):\tanpm\to\tanpm$ is an isomorphism for each $\pnt\in\Fix\flow$.
\end{enumerate}
By definition, put $\Dflow = \Cflow \cap \DiffM$ and let
$\CflowId$ and $\DflowId$ be the identity path components
of the corresponding spaces $\Cflow$ and $\Dflow$ in the $C^{0}_{W}$-topologies.
The next lemma follows from Corollary~\ref{cor:zfix_hdiff} and definitions.
\begin{lem}\label{lm:invshImFn_DflowId}
$\IM\Shift\subset\CflowId$ and $\Shift(\ImFnplus) \subset \DflowId$.\qed
\end{lem}

\begin{defn}
Let $\Diskk$ be the $k$-dimensional disk, $k\in\NNNz$.
A point $\pnt\in\Fix\flow$ is an {\em\condExtk}-point of $\flow$ if
there exists an open neighborhood $\nbs$ of $\pnt$ such that the following holds true:

Suppose $\afunc:(\nbs\setminus\Fix\flow) \times \Diskk \to \RRR$ is a
$C^{\infty}$-function such that the mapping
$\amap:(\nbs\setminus\Fix\flow) \times \Diskk \to \manif$
defined by $\amap(x, t)=\flow(x,\afunc(x,t))$ has a
$C^{\infty}$-extension onto $\nbs\times\Diskk$.
Then $\afunc$ has a $C^{\infty}$-extension onto $\nbs\times\Diskk$.

A point $\pnt$ is {\em\condExt} if it is {\em\condExtk} for any $k \in \NNNz$.
\end{defn}
\begin{thm}\label{th:IMinvsh=CflowId}
Suppose $\Int\Fix\flow=\emptyset$ and each fixed point of $\flow$ is {\em\condExtzero}. Then
\begin{gather}
\label{equ:IMinvsh=CflowId}
\IM\Shift = \CflowId, \\
\label{equ:IMImFn=DflowId}
\Shift(\ImFnplus) = \DflowId.
\end{gather}
\end{thm}
\proof
It is easy to see that~\eqref{equ:IMinvsh=CflowId} implies~\eqref{equ:IMImFn=DflowId}.
So let us prove~\eqref{equ:IMinvsh=CflowId}.
Note that the $C^{0}_{W}$-topologies of the spaces $\smfunc$ and $\smmap$
coincide with the compact-open ones.
This allows us to consider homotopies instead of continuous maps from $I$ into these spaces.

Let $\amap \in \CflowId$.
We will show that there exists a function
$\afunc\in\smfunc$ such that $\amap = \Shift(\afunc)$.
By definition of $\CflowId$, there exists a continuous map $\nu:I \to \smmap$
such that $\nu(0)=\id_{\manif}$ and $\nu(1)=\amap$.
We will construct a map $\tilde\nu: I \to \smfunc$ such that $\nu = \Shift \circ \tilde\nu$.
Then $\amap = \nu(1) = \Shift\circ\tilde\nu(1) \in \IM\Shift$.
Since $\nu(0)=\id_{\manif}$, we can set $\tilde\nu(0)=0$.

In the ``homotopy'' language, the previous paragraph means that
there exists a homotopy $\Amap:\manif \times I \to \manif$
such that $\Amap(x,t) = \nu(t)(x)$, $\Amap(x,0)=x$ and $\Amap(x,1)=\amap(x)$.
Our aim is to construct a homotopy
$\tAfunc:\manif \times I \to \RRR$ such that
$\Amap(x,t) = \Shift(\tilde\nu(t))(x) = \flow(x, \tAfunc(x,t))$ and $\tAfunc(x,0)=0$.
We have that $\tAfunc$ is defined on $\manif\times 0$ by $\tAfunc(x,t)=x$
and we must extend it to $\manif\times I$.

Let $\omega$ be a regular trajectory of $\flow$, $x\in\omega$, and
$p:\RRR \to \omega$ be defined by $p(t)=\flow(x,t)$.
Note that $\Amap(x \times I) \subset \omega$ and $\Amap(x,0)=x$.
Then there exists a {\em unique} function $\tAfunc_{x}:I \to \RRR$
such that $\tAfunc_{x}(0)=0$ and
$F(x,t) = p\circ\tAfunc_{x}(t) = \flow(x, \tAfunc_{x}(t))$ for all $t\in I$.

Indeed, suppose $\omega$ is homeomorphic to the circle $S^1$.
Then $p$ is a covering map and our statements follow from the covering homotopy property for $p$.
If $\omega$ is a non-closed trajectory,
then $p$ is continuous and bijective, though possibly not a homeomorphism.
Nevertheless the restriction of $p$ to any compact subset of $\RRR$ is a homeomorphism,
as being a continuous and injective map from a compact space into a Hausdorff space.
So we put $\tAfunc_{x}=p^{-1}\circ \Amap$.

Define $\tAfunc$ on $(\manif \setminus \Fix\flow) \times I$
by $\tAfunc(x,t) = \tAfunc_{x}(t)$.
Then by formula~\eqref{equ:reg_pnt_lsect}, where $\afunc$ can depend on a parameter,
$\tAfunc(x,t)$ is $C^{\infty}$ in $x$ for each $t\in I$.

Thus for each $t \in I$ the $C^{\infty}$-map $\Amap_{t}$ has
the $C^{\infty}$ shift-function $\tAfunc_{t}$
defined on $\manif\setminus\Fix\flow$.
Let $x\in\Fix\flow$.
Then by the condition \condExtzero\ for $x$,
$\tAfunc_{t}$ can be $C^{\infty}$-extended onto some neighborhood of $x$
remaining a shift-function for $\Amap_{t}$.
Since $\Fix\flow$ is nowhere dense in $\manif$, all these extensions are coherent
and yield a well-defined $C^{\infty}$-extension of $\Amap_{t}$.
In particular $\amap = \Amap_{1} = \Shift(\Afunc_{1})$.
\qed


\subsection{A point that is not {\em\condExtzero}.}
Consider the differential equation on $\RRR$: \ $\frac{d x}{d t}=x^{n}, (n \geq 2)$ and
let $\flow$ be the corresponding local flow defined on the interval $\I = (-a,a)$, $a>0$.
Evidently $\flow$ has exactly three trajectories: $(-a,0)$, $0$ and $(0,a)$.
We will show that the origin $0$ is not a \condExtzero-point of $\flow$.

\proof
Note that the space $\End(\I, \flow)$ consists of $C^{\infty}$-functions $\amap$ on $\I$
preserving the sign of points and such that $\amap'(0)>0$.
Therefore it is path connected in $C^{0}_{W}$-topology, i.e., $\EndId(\I, \flow) = \End(\I, \flow)$.
Let $\amap \in \End(\I, \flow)$.
Then $\amap(0)=0$ and $\amap'(0)>0$.
Therefore, by the Hadamard lemma (see formula~\eqref{equ:Hadamard_proof}),
$\amap(\pnt)=\pnt \bmap(\pnt)$,
where $\bmap$ is a unique $C^{\infty}$-function on $\I$ such that $\bmap(0)=\amap'(0)>0$.

Let us calculate the time $\afunc(\pnt)$ between points $\pnt$ and $\amap(\pnt)$, where $\pnt\in\I$.
\begin{multline*}
\afunc(\pnt) = \int^{\amap(\pnt)}_{\pnt} dt =
\int^{\amap(\pnt)}_{\pnt} \frac{dx}{x^n} =
\frac{ z^{n-1} - \amap(z)^{n-1} }{(n-1) \amap(z)^{n-1} z^{n-1} } =\\
\frac{ 1 - \bmap(z)^{n-1} }{(n-1) \amap(z)^{n-1} } =
\frac{1-\bmap}{\pnt^{n-1}} \cdot \frac{1 + \bmap + \bmap^2 + \cdots + \bmap^{n-2} }{(n-1)\bmap^{n-1}}.
\end{multline*}
It follows that $\afunc$ is $C^{\infty}$ at $0$ if, and only if, $\amap = \pnt + \pnt^{n}\cmap(\pnt)$,
where $\cmap$ is a $C^{\infty}$-function on $\RRR$
(equivalently $\amap(0)=0$, $\amap'(0)=1$ and $\amap^{(k)}(0)=0$ for $k=2,\ldots,n-1$).
Thus for each $n \geq 2$, we have $\IM\Shift \not= \EndId(\I, \flow)$.
\qed

Note that in these cases the flow $\flow$ is not linear.
We will prove in the next section that $\IM\Shift = \EndId(\I, \flow)$ for linear flows.


\section{Regular factors and extensions of flows}\label{sect:calculations}
We prove here that fixed points of ``regular'' extensions of linear flows are \condS\ and \condExt.
Let us represent $\RRR^{m+n}$ as $\RRR^m \times \RRR^n$
and denote its points by $(x,y)$, where $x\in\RRR^m$ and $y\in\RRR^n$.
Let also $\nbm$ and $\nbn$ be open disks in $\RRR^m$ and $\RRR^n$ (respectively)
with centers at the origins and $\nbmn = \nbm \times \nbn$.

\begin{defn}\label{def:reg_fact_ext}
Let $\extflow: \nbmn \times\Interval \to \RRR^{m+n}$
and $\factflow: \nbm \times\Interval \to \RRR^{m}$ be partial flows and
$p_m:\RRR^{m+n} \to \RRR^{m}$ be the natural projection.
Then $\factflow$ is a {\em regular factor} of $\extflow$ and $\extflow$ is,
in turn, a {\em regular extension} of $\factflow$
whenever for each $t \in \Interval$ the following
condition holds:
\begin{equation}\label{equ:reg_ext_diagr}
 p_m \circ \flow_t = \factflow_t \circ p_m .
\end{equation}
\end{defn}
Flows $\extflow$ and $\factflow$ are
{\em regularly equivalent} when they are regular factors of each other.
A flow $\extflow$ is {\em regularly minimal}
if it is nonconstant and each of its regular nonconstant
factors is regularly equivalent to $\extflow$.

Rewriting $\extflow$ in the coordinates $(x,y,t)$
of $\RRR^m \times \RRR^n \times \RRR$ we see
that~\eqref{equ:reg_ext_diagr} is equivalent to the following representation:
\begin{equation}\label{equ:reg_ext_coord}
\extflow(x,y,t) = \left(\factflow(x,t), \B(x,y,t)\right),
\end{equation}
where $\B:\nbmn \times\Interval \to \RRR^n$ is a $C^{\infty}$-map.
Thus $\factflow$ is the ``former'' coordinate function of $\extflow$
and does not depend on $y$.

It follows from the Hadamard lemma that any $C^{\infty}$-map $f:\RRR^n \to \RRR^m$ such that $f(0)=0$
can be represented in the form $f(x) = A(x) \cdot x$, where $x \in \RRR^n$,
and $A(x)$ is a $C^{\infty}$ $(m \times n)$-matrix.
Suppose now that in Definition~\ref{def:reg_fact_ext} the origin $0$
is a fixed point of $\extflow$.
Then formula~\eqref{equ:reg_ext_coord} can also be rewritten in the matrix form as
\begin{equation}\label{equ:reg_ext_matr}
\flow(x, y, t) = \matr{P(x,t)}{0}{Q(x,y,t)}{\Rot{x,y,t}} \cdot \vect{x}{y},
\end{equation}
where $P, Q ,R$ are $C^{\infty}$-matrices of dimensions
$m \times m$, $m \times n$ and $n \times n$, respectively,
such that $P$ does not depend on $y$.
Hence $\factflow(x,t)=P(x,t)x$.

\begin{thm}\label{th:shift_reg_ext_lin}
Let $\extflow:\nbmn \times \Interval \to \RRR^{m+n}$ be a nontrivial partial flow
such that the origin $\pnt_{m+n}=0\in\RRR^{m+n}$ is a fixed point of $\extflow$.
Suppose there exists a linear flow $\factflow(x,t)=e^{At}x$ on $\RRR^m$
such that $\extflow$ is a regular extension of $\factflow$ at $\pnt_{m+n}$.
Then $\pnt_{m+n}$ is {\em\condS} and is an {\em\condExt}-point of $\extflow$.
\end{thm}
\proof
First we show that the properties \condS\ and \condExt\ are inherited by regular extensions
(Lemmas~\ref{lm:condLS_for_extfl} and~\ref{lm:condExt_for_extfl}).
Then we prove them for linear flows.
So let $\extflow$ and $\factflow$ be the flows of Definition~\ref{def:reg_fact_ext},
$\pnt_{m} = 0 \in \RRR^m$ and $\pnt_{m+n} = 0\in\RRR^{m+n}$ be the origins,
$\Diskk$ be an open $k$-disk,
$\nbnk = \nbn\times\Diskk$ and $\nbmnk=\nbm\times\nbn\times\Diskk$.


\begin{lem}\label{lm:condLS_for_extfl}
The origin $\pnt_{m+n}$ is an {\em\condS}-point for $\extflow$
whenever so is $\pnt_{m}$ for $\factflow$.
\end{lem}
\proof
Let $\nbs\subset\nbmnk$ be an open neighborhood of $(\pnt_{m+n}, 0)$
with compact closure $\overline{\nbs}$ and $\afunc \in \smsp{\nbs,\RRR}$.
Then the following diagram is commutative:
$$
\begin{CD}
\smsp{\nbs,\RRR} @= \smsp{\nbs,\RRR} \\
@V{\extflsh}VV @VV{\factflsh}V \\
\smsp{\nbs, \RRR^{m+n}} @>{\tp_m}>> \smsp{\nbs, \RRR^{m}}
\end{CD}
$$
where $\extflsh(\afunc)(x,y,s) = \extflow(x,y, \afunc(x,y,s))$,
$\factflsh(\afunc)(x,y,s) = \factflow(x,\afunc(x,y,s))$, and
$\tp_m(\cmap) = p_m \circ \cmap$
for $\afunc\in\smsp{\nbs,\RRR}$ and $\cmap\in\smsp{\nbs, \RRR^{m+n}}$.

Let $\unifuncext:\nbmn\to(0,\infty)$ and $\unifuncfact:\nbm\to(0,\infty)$ be functions
satisfying the statement of Proposition~\ref{prop:inject_phiV}
for the flows $\extflow$ and $\factflow$ respectively and such that
$\unifuncext(x,y) \leq \unifuncfact(x)$.
Define two $C^{0}_{S}$-neighborhoods of $\afunc$ in the space $\smsp{\nbs,\RRR}$ by
\begin{gather*}
\NBffact = \{ \bfunc\in\smsp{\nbs,\RRR} \ | \
 |\afunc(x,y,s) - \bfunc(x,y,s)| < \unifuncfact(x) \}, \\
\NBfext  = \{ \bfunc\in\smsp{\nbs,\RRR} \ | \
 |\afunc(x,y,s)- \bfunc(x,y,s)| < \unifuncext(x,y) \}.
\end{gather*}
Then $\NBfext\subseteq\NBffact$ and
the restrictions $\factflsh|_{\NBffact}$ and $\extflsh|_{\NBfext}$ are injective.
Hence the inverse mapping $\extflsh^{-1}:\extflsh(\NBfext)\to \NBfext$
coincides with the composition $\factflsh^{-1} \circ \tp_m$.
By the condition \condSnk\ for $\pnt_m$ with respect to $\factflow$, the inverse map
$\factflsh^{-1}:\factflsh(\NBffact) \to \NBffact$ is $C^{\smk}_{W}$-continuous for all $\smk\in\NNNz$.
Since $\tp_m$ is also $C^{\smk}_{W}$-continuous, we see that
so is $\extflsh^{-1}=\factflsh^{-1} \circ \tp_m$.
\qed

\begin{lem}\label{lm:condExt_for_extfl}
The origin $\pnt_{m+n}$ is an {\em\condExt}-point for $\extflow$
whenever so is $\pnt_{m}$ for $\factflow$.
\end{lem}
\proof
Let $\afunc\in\smsp{\nbmnk,\Interval}$ be such that the map
$$
\amap(x,y,s)=\extflow(x,y,\afunc(x,y,s)) =
\left(
       \factflow(x,\afunc(x,y,s)),
       \B(x,y,s)
\right)
$$
has a $C^{\infty}$-extension to $\nbmnk$.
We will show that $\afunc$ has a $C^{\infty}$-extension to $\nbmnk$.
First note that
\begin{equation}\label{equ:rel_between_fix_pnt}
 (\nbm \setminus \Fix\factflow) \times (\nbnk) \ \subset \
 ( \nbmn \setminus \Fix\extflow ) \times \nbk.
\end{equation}
Indeed, it is obvious that $\Fix\extflow \subset \Fix\factflow \times \nbn$.
Then
$$
(\nbm \setminus \Fix\factflow) \times \nbn \ \subset \
  \nbmn \setminus \Fix \extflow.
$$
Multiplying both sides of this relation by $\nbk$ we get~\eqref{equ:rel_between_fix_pnt}.
Since $\pnt_m$ is an \condExtnk-point of $\factflow$, we obtain
that $\afunc$ has a $C^{\infty}$ extension onto $\nbmnk$.
\qed

To complete the theorem
it remains to prove that for each nontrivial linear flow
$\factflow$ on $\RRR^m$ the origin $\pnt_m$ is \condExt\ and \condS.
Notice that we may consider regularly minimal linear flows only.
They are described by the following lemma.
The proof is immediate and will be omitted.
\begin{lem}\label{lm:nontriv_min_reg_fact}
A nonconstant linear flow $\factflow(x,t)=e^{At}x$ is regularly minimal iff
the matrix $A$ is a conjugate to one of the following matrices:
\begin{align*}
& (1) && J_1(\lambda)=\|\lambda\|, \ \lambda\not=0, && \factflow(x,t) = x e^{\lambda t}, \ x\in\RRR; & \\
& (2) && \Rot{\alpha,\beta}=\matr{\alpha}{-\beta}{\beta}{\alpha}, \ \beta \not=0, &&
   \factflow(z,t)=z  e^{(\alpha+\beta i)t}, z\in\CCC; & \\
& (3) && J_2(0)=\matr{0}{1}{0}{0}, && \factflow(x,y,t)=(x+ty, y), \ (x,y)\in\RRR^2. &
\end{align*}
\end{lem}
Let $\Field$ denote either $\RRR$ or $\CCC$,
$\nb$ be an open neighborhood of the origin $0\in\Field$,
$\onb = \nb \setminus \{ 0 \}$ be a ``punctured'' neighborhood of $0$,
$\factflow$ be a linear flow on $\Field$ generated by one of the corresponding
matrices (1)-(3) of Lemma~\ref{lm:nontriv_min_reg_fact}
(in the case (3) we identify $\CCC$ with $\RRR^2$),
and $\factflsh$ be the corresponding shift-map of $\factflow$.

Let $\stfunc: \onb \times \Diskk \to \RRR$ be a $\smsp{}$-function
such that $\stmap(x,\tau) = \factflsh(\stfunc)(x,\tau)=\factflow(x,\tau,\stfunc(x,\tau))$
is $\smsp{}$ on $\nb\times\Diskk$.
We will show that $\stfunc$ can be $\smsp{}$-extended to $\nb\times\Diskk$
and obtain explicit formulas expressing $\stfunc$ in terms of $\stmap$.
They will imply the properties \condS\ and \condExt.
The proof is based on two lemmas.
\begin{lem}\label{lm:HadamardLemmas}
Let $\factflow$ be a flow as in cases (1) or (2) of Lemma~\ref{lm:nontriv_min_reg_fact}.
Then there exists a unique smooth function
$\dmap:\nb\times\Diskk \to \Field$ such that $\stmap(z,\tau) = z \cdot \dmap(z,\tau)$
and $\dmap(0,\tau)\not=0$ for all $\tau\in\Diskk$.
\end{lem}
\begin{lem}\label{lm:ContDivis}
The map
$\mult: \smsp{\nb\times\Diskk, \Field} \to \smsp{\nb\times\Diskk, \Field}$
defined by
$\mult(\stmap)(z,\tau) = z \cdot \stmap(z,\tau)$
is a $C^{\smk}_{W}$-embedding
(i.e., a homeomorphism onto the image in the $C^{\smk}_{W}$-topologies)
for each $\smk\in\NNNz$.
\end{lem}

Assuming that these lemmas hold consider the following cases of $\factflow$.

{\bf (1)}
We have
$\stmap(x,\tau) = \factflow(x, \stfunc(x,\tau)) =
x e^{\lambda \stfunc(x,\tau)}$ for $x\not=0$.
Since $0 \in \Fix\factflow$, it follows that $\stmap(0,\tau)=0$.
From Lemma~\ref{lm:HadamardLemmas} we get $\stmap(x,\tau) = x \dmap(x,\tau)$.
Hence $\dmap(x,\tau) = e^{\lambda \stfunc(x,\tau)}$ and
$\dmap(0,\tau) = \stmap'_x(0,\tau) >0$. Thus
\begin{equation}\label{equ:sh_J1_lambda}
\stfunc(x,\tau) = \frac{1}{\lambda} \ln \dmap(x,\tau).
\end{equation}

{\bf (2)}
Denote $\omega = \alpha + i \beta$.
Then $\stmap(z,\tau)=z e^{\omega\stfunc(z,\tau)}$ for $z\not=0$.
Using Lemma~\ref{lm:HadamardLemmas} we get
$\stmap(z,\tau) = z \dmap(z,\tau)$, where $\dmap$ is $C^{\infty}$.
Hence $\dmap(z,\tau) = e^{\omega \stfunc(z,\tau)}$.
Now we separate the cases of $\alpha$.
If $\alpha \not=0$, then
\begin{equation}\label{equ:sh_Rab_not1}
\stfunc(z,\tau)
= \frac{1}{2\alpha} \ln \left( \dmap(z,\tau) \overline{\dmap(z,\tau)} \right)
= \frac{1}{2\alpha} \ln \left| \dmap(z,\tau) \right|^2.
\end{equation}
Suppose $\alpha=0$.
Then $\dmap(z,\tau) = e^{i \beta \stfunc(z,\tau)}$.
It follows that $\stfunc$ is not unique and is determined up to a constant summand,
\begin{equation}\label{equ:sh_Rab_1}
\stfunc(z,\tau) = \frac{1}{\beta} \arg ( \dmap(z,\tau) ) +
 \frac{2\pi k}{\beta}, \ k \in \ZZZ.
\end{equation}

{\bf (3)}
In this case, $\Fix\factflow = \{(x,0)\ | \ x\in\RRR\}$
and $$\stmap(x,y,\tau) = \factflow(x,y,\stfunc(x,y,\tau)) =
 (x+ y \stfunc(x,y,\tau), y)$$ for $y\not=0$.
Let $\stmap_1(x,y,t)=x+y\stfunc(x,y,\tau)$ be the first coordinate
function of $\stmap$.
Then $\stfunc = (\stmap_1-x)/y$ for $y\not=0$.
Note that the function $\Bmap(x,y,\tau) = \stmap_1(x,y,\tau) - x$
is $C^{\infty}$ and $\Bmap(x,0,\tau)\equiv 0$.
Therefore we can apply the Hadamard lemma to $\Bmap$ and obtain
a unique $C^{\infty}$ function $\dmap: \nb\times\Diskk \to \RRR$
such that $\Bmap(x,y,\tau)=y\dmap(x,y,\tau)$. Hence
\begin{align}\label{equ:sh_J2_0}
\stfunc(x,y,\tau) = \frac{\stmap_2(x,y,\tau) - x}{y} = \dmap(x,y,\tau).
\end{align}

From formulas~\eqref{equ:sh_J1_lambda}-\eqref{equ:sh_J2_0}
we see that in all cases the function $\stfunc$
has a $C^{\infty}$ extension to some neighborhood of $0\in\Field$.
This implies the condition \condExt\ for $\factflow$ at $0\in\Field$.
Furthermore, let $\nbs$ be any small open neighborhood $\nbs$ of $0$.
It follows from Lemma~\ref{lm:ContDivis} and these formulas
that there exists a $C^{0}_{W}$-neighborhood of $\stmap|_{\nbs}$
in $\smsp{\nbs,\RRR}$ such that the correspondence
$\stmap|_{\nbs} \,\mapsto\, \dmap|_{\nbs} \,\mapsto\,\stfunc|_{\nbs}$
is $C^{\smk}_{W}$-continuous.
Hence $0$ is an $\condS$-point of $\factflow$.
To complete the proof Theorem~\ref{th:shift_reg_ext_lin}, we must prove
Lemmas~\ref{lm:HadamardLemmas} and~\ref{lm:ContDivis}.

\proofstyle{Proof of Lemma~\ref{lm:ContDivis}.}
Note that $\mult$ is linear, injective, and $C^{\smk}_{W}$-continuous
for any $\smk\in\NNNz$.
Let us verify the $C^{\smk}_{W}$-continuity of the inverse map $\mult^{-1}$.

For each compact set $K \subset \nb \times \Diskk$ and
$\smk\in\NNNz$, consider the norm
$\normnK{\cdot}$ on $\smsp{\nb\times\Diskk, \Field}$ defined by
$\normnK{\stmap} =
 \sum\limits_{i=0}^{\smk} \sup\limits_{x \in K} {|D^{i}\stmap(x)|}$,
where $|D^{i}\stmap(x)|$ denotes the sum of absolute values
of all derivatives of $\stmap$ of degree $i$.
If $K$ runs over all compact subsets of $\nb \times \Diskk$,
the norms $\normnK{\cdot}$ generate the $C^{\smk}_{W}$-topology
in $\smsp{\nb\times\Diskk, \Field}$.
Let $\stmap \in \smsp{\nb\times\Diskk, \Field}$ and
$\dmap = \mult(\stmap) = z \stmap$.
The following inequality implies our lemma and can be easily verified:
\begin{multline*}
\normnK{\dmap} \ \leq \ |z| \normnprevK{\stmap} + \normnK{\stmap} \ \leq \\
(1 + |z|) \normnK{\stmap} \ \leq \ (1+\diam K) \normnK{\stmap}. \qed
\end{multline*}


\section{Proof of Lemma~\ref{lm:HadamardLemmas}}\label{sect:proof_HadamLemmas}
{\bf (1)}
In this case the lemma follows from the well-known Hadamard lemma.
Indeed, since $\stmap(x,\tau) = e^{\lambda \stfunc(x,\tau)}x$ is $C^{\infty}$,
we see that $\stmap(0,\tau)=0$ for all $\tau$.
Then
\begin{equation}\label{equ:Hadamard_proof}
\stmap(x,\tau) =
  \int\limits_{0}^{x} \frac{\partial\stmap}{\partial t}(t,\tau) dt =
  x \int\limits_{0}^{1} \frac{\partial\stmap}{\partial t}(t \cdot x,\tau) dt.
\end{equation}
Denoting the last integral by $\dmap(x,\tau)$, we get
that $\stmap(x,\tau)=x\dmap(x,\tau)$,
$\dmap$ is $C^{\infty}$,
and $\dmap(0,\tau) = \stmap'(0,\tau) \not=0$ for each $\tau\in\Diskk$.

{\bf (2)}
Let us denote $\omega = \alpha + i \beta$.
Then $\stmap(z,\tau) = e^{\omega\stfunc(z,\tau)}z$, so we must put
\begin{equation}\label{equ:phi_by_shift}
\dmap(z,\tau) = e^{\omega\stfunc(z,\tau)}, \forall z\not=0.
\end{equation}
Hence
\begin{equation}\label{equ:sh_atan_ln}
\stfunc = \frac{1}{2 \alpha}\ln |\dmap|^2 = \frac{1}{\beta}\arg \dmap,
\quad \forall z\not=0.
\end{equation}
\begin{lem}
The function $\dmap$ satisfies the following equation:
\begin{equation}\label{equ:rel_for_phi}
\IM(\omega \, \dmap \, d\bar\dmap) = 0.
\end{equation}
\end{lem}
\proof
From~\eqref{equ:phi_by_shift} we get $d\dmap = \omega \dmap d\stfunc.$
Multiplying both sides of this formula by $d\bar\dmap$
and taking into account that $d\stfunc$ and $d\dmap d\bar\dmap$ are real,
we see that so is $\omega \, \dmap \, d\bar\dmap$.
\qed

To complete the proof of our lemma we separate the cases $\alpha\not=0$ and $\alpha=0$.
\begin{lem}
If $\alpha \not= 0$, then the functions $\stfunc$ and $\dmap$ are $C^{\infty}$.
\end{lem}
\proof
It suffices to prove that $|\dmap(z)|^2$ is $C^{\infty}$.
Indeed, since $\stmap$ is a diffeomorphism at $0$,
there exist constants $c$ and $C$ such that $0 < c < |\stmap(z)|/|z| = |\dmap(z)| < C$
in some neighborhood of $0\in\CCC$.
Thus, if $|\dmap|^2$ is $C^{\infty}$, then by~\eqref{equ:sh_atan_ln}
and~\eqref{equ:phi_by_shift} so are $\stfunc$ and $\dmap$.

Now, let us expand formula~\eqref{equ:rel_for_phi},
$$
\omega \, \dmap \, d\bar\dmap =
\omega \ \frac{\stmap}{z} \ d \left( \frac{\bar\stmap}{\barz} \right) =
\frac{ \omega \, \stmap}{z} \cdot \frac{\barz \, d\bar\stmap - \bar\stmap\,d\barz}{\barz^2} =
\frac{ z\barz \cdot \omega \, \stmap \, d\bar\stmap - \stmap \bar\stmap \cdot \omega \, z \, d\barz }
  {(z \barz)^2}\, .
$$
Since $z\barz$, $\stmap\bar\stmap$ are real,
the relation~\eqref{equ:rel_for_phi} is equivalent to the following one:
$$
\stmap \bar\stmap \cdot \IM( \omega \, z \, d\barz ) =
 z\barz \cdot \IM( \omega \, \stmap \, d\bar\stmap ),
$$
whence
\begin{equation}\label{equ:g2imz_imh}
|\dmap|^2 \cdot  \im( \omega \, z \, d\barz ) = \im( \omega \, \stmap \, d\bar\stmap ).
\end{equation}
Substituting $d\barz = \baromega = \alpha - i\beta$
in the last formula we get
$\im( \omega \, z \, \baromega ) = y |\omega|^2$.
Then the left side of equation~\eqref{equ:g2imz_imh} becomes equal to
$ |\dmap(x,y)|^2 \cdot y |\omega|^2$.
This function is $\smsp{}$ and so is the right-hand side.
It follows from the Hadamard lemma that $|\dmap|^2$ is also $\smsp{}$.
\qed


Suppose now that $\alpha=0$.
Then $\factflow(z,t) = e^{i \beta t} z$ and $\beta \not=0$.
Therefore,
\begin{equation}\label{equ:p2+q2=x2+y2}
z \barz = \stmap \bar\stmap.
\end{equation}
In fact, this is just another expression of~\eqref{equ:rel_for_phi} for our case $\alpha=0$.
\begin{claim}\label{clm:damap_dbarz0_0}
$
\displaystyle{
        \frac{\partial^n\!\stmap}{\partial \bar z^{n}}(0) = 0
}
$ for each $n=1,2,\ldots$
\end{claim}
For the proof we need the following lemma.
\begin{lem}\label{lm:f=f_dh0}
Let $\stmap:\RRR^n \to \RRR^n$ be a $C^1$-diffeomorphism
such that $\stmap(0)=0$.
Let $\stmap'(0) : \RRR^n \to \RRR^n$ be the tangent map of $\stmap$ at $0$,
and $\bmap:\RRR^n \to \RRR$ be a continuous
homogeneous function of degree $k$, i.e.,
$\bmap(tx)=t^k \bmap(x)$ for $t>0$ and $x \in \RRR^n$.
If $\bmap = \bmap \circ \stmap$ then $\bmap = \bmap \circ \stmap'(0)$.
\end{lem}
\proof
Let $x \in \RRR^n$ and $t>0$. Then
$$\bmap(x)=\frac{\bmap(t  x)}{t^k} = \frac{\bmap(\stmap(t  x))}{t^k} =
\bmap \left( \frac{\stmap(t  x)}{t} \right) \
\mathop\to_{t \to 0} \ \bmap \circ \stmap'(0)(x). \ \qed$$

\proofstyle{Proof of Claim~\ref{clm:damap_dbarz0_0}.}
Let $\stmap(z) = p(z) + i q(z)$, where $p, q \in \smsp{\CCC,\RRR}$.
We will use the induction in $n$.
Let $n=1$.
Then the relation \eqref{equ:p2+q2=x2+y2} means that
$\stmap$ preserves the homogeneous polynomial $\bfunc(x,y)=x^2 + y^2$.
It follows from Lemma~\ref{lm:f=f_dh0} that so does $\stmap'(0)$.
Hence $\stmap'(0)$ is an orthogonal matrix and $\frac{\partial \stmap}{\partial \bar z}(0)=0$.
Thus $\stmap'_z(0)$ coincides with multiplication by $e^{i a}$, where $a \in [0, 2\pi)$,
whence $\stmap(z) = e^{i a} z + \varepsilon_{2}$, where $\varepsilon_{2} = o(|z|^2)$.
Let us define $\dmap(0)=e^{i a}$. Then $\dmap$ becomes continuous at $0$.

Suppose we have proved the lemma for $n-1$.
Then
$$ \stmap(z) = e^{ia} z + A(z) z + b \bar z^{n} + \varepsilon_{n+1},$$
where $A(z)$ is a polynomial in the variables $z$ and $\bar z$,
$1 \leq \deg A \leq n-1$, and
$
  b = \frac{1}{n!} \, \frac{\partial^n \stmap}{\partial \bar z^n}(0).
$
Thus we have to prove that $b=0$.

Substituting $h$ in formula~\eqref{equ:p2+q2=x2+y2} we obtain
$$ z \barz = \stmap \bar\stmap = z \barz + A(z) z \barz + b \bar z^{n+1} + \bar A(z) z \bar z + \bar b  z^{n+1} + \theta_{n+2}, $$
where $\theta_{n+2} = o(|z|^{n+2})$.
Hence
$$
 \bar b z^{n+1}+ (A(z) + \bar A(z)) z \barz + b \bar z^{n+1}= 
 -\theta_{n+2} 
$$
The right part of this equality is a function of order $|z|^{n+2}$,
while the left one is a polynomial of degree $\leq n+1$ in the variables $z$ and $\barz$.
Hence all the coefficients of this polynomial are zeros.
Since the first two summands contain a multiple $z$,
we obtain that the coefficient at $\bar z^n$ is $b$.
Hence $b=0$. \qed

The following lemma is left to the reader.
\begin{lem}\label{lm:smooth_e/z}
Let $\varepsilon: \CCC \to \CCC$ be a $C^{\infty}$-function such that $\varepsilon = o(|z|^k)$.
Define $\tau:\CCC \to \CCC$ by $\tau(z) = \varepsilon(z)/z$ for $z\not=0$ and $\tau(0)=0$.
Then $\tau$ is $C^{k-2}$. \qed
\end{lem}
Now we can complete the case (2).
It follows from Claim~\ref{clm:damap_dbarz0_0} that
for each $n \in \NNN$ the Taylor expansion of degree $n$ of $\stmap$ at $0$ has the form
$\stmap(z) = \teyl_{n-1} z + \varepsilon_{n+1},$
where $\teyl_{n-1}$ is a polynomial of degree $n-1$ in the variables $z$ and $\bar z$
and $\varepsilon_{n+1} = o(|z|^{n+1})$.
Hence $\dmap(z) = \stmap(z)/z = \teyl_{n-1} + \varepsilon_{n+1}/z.$
Applying Lemma~\ref{lm:smooth_e/z} to the remainder $\varepsilon_{n+1}/z$,
we see that this function is $C^{n-1}$.
Hence so is $\dmap$ for any $n\in\NNN$, i.e., $\dmap$ is $C^{\infty}$.
Lemma~\ref{lm:HadamardLemmas} is proved.
\qed

\section{Proof of \MainTheorem}\label{sect:HomotopyGroups}
Let $\flow$ be a global flow on $\manif$
such that for each fixed point $\pnt$
of $\flow$ there exist local coordinates $(x_1,\ldots,x_n)$ and a nontrivial linear
flow $\factflow$ on $\RRR^{m} (0< m \leq n)$ such that $\pnt=0$ and for all $t$
in some neighborhood of $0\in\RRR$ we have
$p_m \circ \flow_t = \factflow_t \circ p_m$,
where $p_m:\RRR^n \to \RRR^m$ is a natural projection.

Then $\Int\Fix\flow=\emptyset$, since this holds for linear flows.
Moreover, it follows from
Theorem~\ref{th:shift_reg_ext_lin} that $\flow$ satisfies
the conditions \condExt\ and \condS\ at each point $\pnt\in\Fix\flow$.
Then by Theorems~\ref{th:LocSect} and~\ref{th:IMinvsh=CflowId}
we obtain that $\IM\Shift = \CflowId$, the map
$\Shift:\smfunc\to\IM\Shift$ is covering, and its group of covering
slices is either $0$ or $\ZZZ$.
Finally, by Theorem~\ref{th:IMinvsh=CflowId} and
Lemma~\ref{lm:ImFnplus_minus_convex},
the set $\ImFnplus = \Shift^{-1}(\DflowId)$ is a convex
subset of the linear space $\smfunc$.

Suppose now that $\manif$ is compact.
Let $X$ denote either $\smfunc$ or $\ImFnplus$.
Then the image $Y=\Shift(X)$ is either $\CflowId$ or $\DflowId$, respectively.
Evidently $X$ is a \Frechet\ manifold, whence so is $Y$, as being its image
under the covering map $\Shift$.
It follows that $Y$ has the homotopy type of a CW-complex, (e.g., Palais~\cite{Palais_RS}).
Since $X$ is contractible,
we get that $Y$ is aspherical, i.e., $\pi_n(Y)=0$ for all $n\geq 2$.
By Theorem~\ref{th:Zid_descr}, $\pi_1(Y)$ is either $0$ or $\ZZZ$.
Thus $Y$ is either contractible or homotopically equivalent to $S^1$.

Since $\Zid = \Shift^{-1}(\id_{\manif}) \subset \ImFnplus$,
we see that the embedding $\DflowId \subset \CflowId$ induces an
isomorphism of all homotopy groups and is, therefore, a homotopy equivalence.
Finally, suppose that either $\flow$ has at least one non-closed trajectory,
or the tangent linear flow at some fixed point of $\flow$ is trivial.
Then by Proposition~\ref{pr:suffcond_Zid_zero}, $\Zid=0$, whence
$\DflowId$ and $\CflowId$ are contractible.
\qed

\section{The closure of $\CflowId$}\label{sect:CflowId_closure}
By \MainTheorem, in most cases the sets $\CflowId$ and $\DflowId$ are contractible.
Nevertheless, as the referee of this paper noted, their closures are likely not so.
For example, in the article~\cite{Keesling} of J.~Keesling
the homeomorphisms group $G$ of a solenoid $\Sigma$ is considered.
The path-component $C$ of the unit element $e$ in $\Sigma$ is
a dense one-parametric subgroup $\phi:\RRR\to\Sigma$ such that $\phi$ is one-to-one.
Then $C$ and $\overline{C} = \Sigma$ are of different homotopy types.
It is proven that the identity path component $G_{\id}$ of $G$
is homotopically equivalent to $C$ while the closure of $G_{\id}$ has the homotopy type of $\Sigma$.

Notice that for each closed subset $K\subset\manif$ the set
$$
 \Kinv(K)=\{ \amap \in \smmap \ | \ \amap(K)\subset K \}
$$
is closed in $\smmap$ with the topology of point-wise convergence.
Therefore it is closed in each Whitney topology of $\smmap$.
Thus, if $\flow$ is a global flow on $\manif$, then the set
$\FlowInv = \mathop\cap\limits_{\omega} \Kinv(\overline{\omega})$,
where $\omega$ runs over all trajectories of $\flow$,
is closed in the Whitney topologies of $\smmap$.
Clearly, $\Cflow\subset\FlowInv$.
Hence $\CflowId\subset\FlowInvId$, where $\FlowInvId$ is the identity path
component of $\FlowInv$ in compact-open topology.
The following lemma is not hard to prove.
\begin{lem}
Let $\flow$ be a {\em global\/} flow defined on
a {\em connected} subset of $\RRR$.
Then $\overline{\CflowId} = \FlowInvId$.
\end{lem}
However, in general, it seems to be a problem to prove this equality
as well as to establish that $\FlowInvId$ is closed.
Consider, for instance, an irrational flow $\flow$ on the $n$-torus $T^n$.
Each trajectory of $\flow$ is everywhere dense in $T^n$, whence $\FlowInv=\smsp{T^n, T^n}$.
Let $d$ be a metric on $\smsp{T^n, T^n}$ yielding the compact-open topology.
Since $T^n$ is ANR, it follows that two continuous mappings $\amap,\bmap:T^n\to T^n$
are homotopic provided $d(\amap,\bmap)$ is sufficiently small.
Therefore each path-component of $\smsp{T^n, T^n}$ is open.
Hence it is also closed as the complement to the union of all other ones.
Thus $\FlowInvId$ is closed.

Notice that $\flow$ has no fixed points.
Therefore it satisfies the conditions of \MainTheorem, whence $\IM\Shift = \CflowId$.
Thus the statement $\overline{\CflowId}=\FlowInvId$ would mean that
$\overline{\IM\Shift}$ is the unity path-component of $\smsp{T^n,T^n}$,
i.e., that for any $\eps>0$ each smooth mapping $\amap:T^n\to T^n$ that is homotopic to $\id_{T^n}$
could be $\eps$-approxima\-ted in metric $d$
by a map of the form $\amap_{\eps}(x)=\flow(x,\afunc_{\eps}(x))$,
where $\afunc_{\eps}\in\smsp{T^n,\RRR}$.
The author does not know whether this is true or not.
One of the difficulties is exposed by the following general proposition:
if $\amap(\pnt)$ does not belong to the trajectory
of $\pnt$, then roughly speaking, $\lim\limits_{\eps\to 0}\afunc_{\eps}(\pnt) = \infty$.

\begin{prop}
Let $\flow$ be a global flow on $\manif$.
Suppose that there exists a trajectory $\omega$ of $\flow$ such that
$\Int\,\overline{\omega} \not= \emptyset$.
Let also $\amap \in \overline{\IM\Shift}\setminus\IM\Shift$, and $\pnt\in\Int\,\overline{\omega}$
be such that $\amap(\pnt)$ does not belong to the trajectory $\omega_{\pnt}$
of $\pnt$.
If $\{\tm_i\}_{i\in\NNN}$ is a sequence of reals such that
$\lim\limits_{i\to\infty}\flow(\pnt,t_i) = \amap(\pnt)$,
then $\lim\limits_{i\to\infty}\tm_i = \infty$.
\end{prop}
\proof
First note that $\omega_{\pnt}$ is non-closed.
Denote $y=\amap(\pnt)$ and $y_i=\flow(\pnt,\tm_i)$ for all $i\in\NNN$.
Fix any $\numb>0$ and define the compact subset
$\omega_{\numb} \subset \omega_{\pnt}$
by $\omega_{\numb} = \flow(\pnt\times[-\numb,\numb])$.
Then $y\not\in\omega_{\numb}$.
Hence there exists a neighborhood $\nb$ of $y$ such that
$\nb\cap\omega_{\numb}=\emptyset$.
Since $\lim\limits_{i\to\infty}y_i = y$ and $y_i\in\omega$,
we have $y_i=\flow(\pnt,\tm_i)\in\omega_{\pnt}\setminus\omega_{\numb}$
whence $\tm_i>\numb$ for almost all $i$.
Taking the number $\numb$ arbitrary large, we obtain that $\lim\limits_{i\to\infty}\tm_i = \infty$.
\qed

\newcommand\Bibitem[5]{
\bibitem[#1]{#2}
{\sc #3}
\emph{#4}
#5
}

\end{document}